\documentclass[aop,MSNbibl,nameyear,dvips]{arximspdf}

\doi{10.1214/11-AOP711} 
\volume{41}
\issue{2}
\pubyear{2013}
\firstpage{785}
\lastpage{816}

\makeatletter

\newproclaim{mydef}{Definition}

\renewcommand{\Pr}{\operatorname{Pr}}

\newtheorem{theorem}{Theorem}
\newtheorem{lem}{Lemma}
\newtheorem{cor}{Corollary}
\newtheorem{prop}{Proposition}

\newtheorem{defin}[mydef]{Definition}
\newproclaim{rem}{Remark}

\newcommand{\lra}{\longrightarrow}
\newcommand{\diam}{\operatorname{diam}}
\newcommand{\Ex}{\mathbb{E}}
\newcommand{\card}{\operatorname{Card}}
\newcommand{\E}{{\mathbb E}}

\makeatother

\begin{document}
\begin{frontmatter}

\title{A CLT for empirical processes involving time-dependent data}
\runtitle{Time-dependent empirical processes}

\begin{aug}
\author[A]{\fnms{James} \snm{Kuelbs}\ead[label=e1]{kuelbs@math.wisc.edu}},
\author[A]{\fnms{Thomas} \snm{Kurtz}\thanksref{t2}\ead[label=e2]{kurtz@math.wisc.edu}}
\and
\author[B]{\fnms{Joel} \snm{Zinn}\corref{}\thanksref{t3}\ead[label=e3]{jzinn@math.tamu.edu}}
\runauthor{J. Kuelbs, T. Kurtz and J. Zinn}
\affiliation{University of Wisconsin, University of Wisconsin and Texas~A\&M~University}
\address[A]{J. Kuelbs\\
T. Kurtz\\
Department of Mathematics\\
University of Wisconsin \\
Madison, Wisconsin 53706\\
USA\\
\printead{e1}\\
\hphantom{E-mail: }\printead*{e2}}
\address[B]{J. Zinn\\
Department of Mathematics \\
Texas A\&M University \\
College Station, Texas 77840\\
USA\\
\printead{e3}} 
\end{aug}

\thankstext{t2}{Supported in part by NSF Grant DMS-08-05793.}

\thankstext{t3}{Supported in part by NSA Grant H98230-08-1-0101.}

\received{\smonth{8} \syear{2010}}
\revised{\smonth{6} \syear{2011}}

%
\begin{abstract}
For stochastic processes $\{X_t\dvtx t \in E\}$, we establish sufficient
conditions for the empirical process based on $\{ I_{X_t \le y} -
\Pr(X_t \le y)\dvtx t \in E, y \in\mathbb{R}\}$ to satisfy the CLT
uniformly in $ t \in E, y \in\mathbb{R}$. Corollaries of our main
result include examples of classical processes where the CLT holds, and
we also show that it fails for Brownian motion tied down at zero and
$E= [0,1]$.
\end{abstract}

%

\end{frontmatter}

\section{Introduction}\label{sec1}

To form the classical empirical process, one starts with i.i.d. random
variables $\{X_{j}\dvtx j\ge1\}$, with distribution function $F$, and with
%
%
\begin{equation}\label{emp-proc}
\mathbb{P}_{n}(A)=\frac1{n}\sum_{j=1}^{n}I_{X_{j}\in A},
\end{equation}
one considers the process $\mathbb{F}_{n}(y)=\mathbb{P}_{n}((-\infty,y])$.

By the classical Glivenko--Cantelli theorem,
\[
\sup_{y\in\mathbb{R}}|\mathbb{F}_{n}((-\infty,y])-F(y)|\lra0\qquad
\mbox{a.s.}
\]
By Donsker's theorem,
\[
\bigl\{\sqrt{n}\bigl(\mathbb{F}_{n}(y)-F(y)\bigr)\dvtx y\in\mathbb{R}\bigr\}
\]
converges in distribution in a sense described more completely below.
Hence limit theorems for such processes, such as the law of large
numbers and the central limit theorem (CLT), allow one to
asymptotically get uniform estimates for the unknown cdf, $F(y)=\Pr
(X\le y)$, via the sample data.

A more general version of these processes is to replace the indicators
of half-lines in (\ref{emp-proc}) by functions of a ``random variable''
taking values in some abstract space $(S,\mathcal{S})$. More specifically,
%
%
\begin{equation}
\Biggl\{\frac{1}{\sqrt{n}}\sum_{j=1}^n \bigl(f(X_j)-\E f(X)
\bigr)\dvtx f\in\mathcal{F}\Biggr\},
\end{equation}
where the index set, $\mathcal{F}$, is a subset of $\mathcal
{L}_{\infty
}(S,\mathcal{S})$ or an appropriate subset of $\mathcal
{L}_{2}(S,\mathcal{S},P)$. We use the notation $\mathcal{L}_p(S,
\mathcal{S},P), 0< p < \infty$, to denote the $\mathcal{S}$-measurable
functions on $S$ whose absolute value to the $p$th power is integrable
with respect to $P$, rather than the equivalence classes of these
functions. Of course, when $p=\infty$ the functions are $\mathcal
{S}$-measurable and uniformly bounded on $S$. The standard notation
$L_p(S, \mathcal{S},P)$ is used when we are dealing with equivalence classes.

However, even in the case that the class of functions is a class of
indicators, unlike the classical case, it is easy to see there are many
classes of sets, $\mathcal{C}$, for which the limit theorem does not
hold. As a matter of fact, the limiting Gaussian may not be continuous,
for example, if $\mathcal{C}=$ all Borel sets of $\mathbb{R}$ or even
$\mathcal{C}=$ all finite sets of $\mathbb{R}$. And further, even if
the limiting Gaussian process is continuous, the limit theorem may
still fail.

Luckily, in the case of sets, modulo questions of measurability, there
are necessary and sufficient conditions for this sequence to converge
in distribution to some mean zero Gaussian process. However, all the
nasc's are described in terms of the asymptotic behavior of a
complicated function of the sample, $\{X_{n}\}_{n=1}^{\infty}$. What we
attempt to do in this paper is to obtain additional sufficient
conditions that are useful when $X$ takes values in some function space
$S$, and the sets in $\mathcal{C}$ involve the time evolution of the
stochastic process $X$. Of course, $\mathcal{C}$ is still a class of
sets, but a primary goal that emerges here is to provide sufficient
conditions for a uniform CLT in terms of the process $X=\{X(t)\dvtx t
\in E\}$ that depend as little as possible on the parameter set $E$.
However, classes of sets such as this rarely satisfy the Vapnik--\v
{C}ervonenkis condition. Also, this class of examples arises naturally
from the study of the median process for independent Brownian Motions
[see Swanson (\citeyear{swanson-scaled-median,swanson-fluctuations})],
where he studies the limiting quantile process for independent Brownian
motions. This was observed by Tom Kurtz, and the follow-up questions
led us to start this study. Here we concentrate on empirical process
CLTs, and our main result is Theorem~\ref{mainresult} below. Another
theorem and some examples showing the applicability of Theorem
\ref{mainresult} appear in Section~\ref{sec7}. In Section~\ref{sec8}
there are additional examples which show some obvious conjectures one
might be tempted
to make, concerning the CLT formulated here, are false.
In particular, the examples in Section~\ref{sec84} motivate the various
assumptions we employ in Theorem~\ref{mainresult}. As for future work,
an upgraded version of \citet{vervaat-quantiles} would perhaps
allow one to relate the results obtained here and those of Swanson, but
this is something to be done elsewhere.

\section{Previous results and some definitions}\label{sec2}

Let $(S, \mathcal{S},P)$ be a probability space, and define $(\Omega,
\Sigma,Pr)$ to be the infinite product probability space
$(S^N,\mathcal
{S}^N,P^N)$. If $X_j\dvtx\Omega\rightarrow S$ are the natural projections
of $\Omega$ into the $j$th copy of~$S$, and $\mathcal{F}$ is a
subset of $\mathcal{L}^2(S,\mathcal{S},P)$ with
%
%
\begin{equation}
{\sup_{ f \in\mathcal{F}}} |f(s)|< \infty,\qquad s \in S,
\end{equation}
then we define
%
%
\begin{equation}
\nu_n(f)=\frac{1}{\sqrt{n}}\sum_{j=1}^n \bigl(f(X_j)-\E
f(X)\bigr),\qquad f \in\mathcal{F}.
\end{equation}
Let $\ell_{\infty}(\mathcal{F})$ be the bounded real valued functions
on $\mathcal{F}$, with the $\sup$-norm, and recall that a Radon measure
$\mu$ on $\ell_{\infty}(\mathcal{F})$ is a finite Borel measure which
is inner regular from below by compact sets. Then the functions $f
\rightarrow f(X_j)- \Ex(f(X_j))$, \mbox{$j \geq1$}, are in $\ell_{\infty
}(\mathcal{F})$, and we say $\mathcal{F} \in \operatorname{CLT}(P)$ if the stochastic
processes $\{\nu_n(f)$, $f \in\mathcal{F}\}, n \ge1$, converge weakly to
a centered Radon Gaussian measure $\gamma_P$ on $\ell_{\infty
}(\mathcal
{F})$. More precisely, we have the following definition.
%
%
\begin{defin}\label{defin1} Let $\mathcal{F}\subset\mathcal{L}_{2}(P)$ and satisfy
(3). Then $\mathcal{F} \in \operatorname{CLT}(P)$, or $\mathcal{F}$ is a $P$-Donsker
class if there exists a centered Radon Gaussian measure $\gamma_P$ on
$\ell_{\infty}(\mathcal{F})$ such that for every bounded continuous
real valued function $H$ on $\ell_{\infty}(\mathcal{F})$, we have
\[
\lim_{n \rightarrow\infty} \Ex^{*}(H(\nu_n))= \int H\,d\gamma_P,
\]
where $\Ex^{*}H$ is the usual upper integral of $H$. If $\mathcal{C}$
is a collection of subsets from $\mathcal{S}$, then we say $\mathcal{C}
\in \operatorname{CLT}(P)$ if the corresponding indicator functions are a $P$-Donsker class.
\end{defin}

The probability measure $\gamma_P$ of Definition~\ref{defin1} is obviously the law
of the centered Gaussian process $G_P$, indexed by $\mathcal{F}$ having
covariance function
\[
\Ex G_{P}(f)G_{P}(g)=\Ex_{P}f g-\Ex_{P}f\Ex_{P}g,\qquad f ,g \in\mathcal{F},
\]
and $L^2$ distance
\[
\rho_{P}(f,g)=\Ex_{P}\bigl(\{(f-g)-\Ex_{P}(f-g)\}^2\bigr)^{{1}/{2}},\qquad f,g
\in
\mathcal{F}.
\]
Moreover, if $\gamma_P$ is as in Definition~\ref{defin1}, then it is known that
the process $G_P$ admits a version all of whose trajectories are
bounded and uniformly $\rho_P$ continuous on $\mathcal{F}$. Hence we
also make the following definition.
%
%
\begin{defin} A class of functions $\mathcal{F}\subset\mathcal
{L}_{2}(P)$ is said to be $P$-pre-Gaussian if the mean zero Gaussian\vadjust{\goodbreak}
process $\{G_{P}(f)\dvtx f\in\mathcal{F}\}$ with covariance and $L_2$
distance as indicated above
has a version with all the sample functions bounded and uniformly
continuous on $\mathcal{F}$ with respect to the $L_2$ distance $\rho_P(f,g)$.
\end{defin}

Now we state some results which are useful for what we prove in this
paper. The first is important as it helps us establish counterexamples
to natural conjectures one might make in connection to our main result,
appearing in Theorem~\ref{mainresult} below.
%
%
\begin{theorem}[{[\citet{empirical-specialinvited} for
sufficiency and \citet{talagrand-donsker-sets} for necessity of the
$\Delta^{\mathcal{C}}$ condition in (ii)]}]
\label{delta-c}
Let $\Delta^{\mathcal{C}}(A)$
denote the number of distinct subsets of $A$ obtained when one
intersects all sets in $\mathcal{C}$ with $A$. Then, modulo
measurability assumptions, conditions \textup{(i)} and
\textup{(ii)} below are equivalent.
\begin{longlist}
\item The central limit theorem holds for the process
\[
\Biggl\{\frac1{\sqrt{n}}\sum_{j=1}^n [I_{X_j\in C} - \Pr(X\in
C)
]\dvtx C\in\mathcal{C}\Biggr\}
\]
or more briefly $\mathcal{C}\in \operatorname{CLT}(P)$.
\item
%
%
\begin{eqnarray}\label{delta-condition}
&&\mbox{\textup{(a)}}\quad \frac{\ln\Delta^{\mathcal{C}}(\{X_1,\ldots,X_n\})}{\sqrt
{n}}\to
0 \qquad\mbox{in (outer) probability}\quad \mbox{and} \\
&&\mbox{\textup{(b)}}\quad
\mathcal{C}\mbox{ is } P\mbox{-pre-Gaussian}.
\end{eqnarray}
\end{longlist}
\end{theorem}

A sufficient condition for the empirical CLT, which is used in the
proof of our main theorem, is given in Theorem 4.4 of \citet{clt-localcondit}.
%
%
\begin{theorem}[{[\citet{clt-localcondit}]}]\label{agoz} Let
\[
\mathcal{F}\subset
\mathcal{L}_{2}(S,\mathcal{S},P),\qquad F={\sup_{f\in\mathcal
{F}}}|f(X)|
\]
and $P$ be the distribution of $X$ with respect to $\Pr$, that is,
$P={\Pr} \circ {X^{-1}}$. Also, let $Pf = \int f(x) P(dx)$. Assume that $
\| Pf\| _{\mathcal{F}}\equiv{\sup_{ f \in\mathcal{F}}} |P(f)| <\infty$
and:
\begin{longlist}
\item $u^{2}\Pr^{*}(F >u)\to0 \mbox{ as } u\to\infty$;
\item $\mathcal{F}$ is $P$-pre-Gaussian;
\item there exists a centered Gaussian process $\{G(f)\dvtx f \in
\mathcal{F}\}$ with $L_2$ distance $d_G$ such that $G$ is sample bounded
and uniformly $d_G$ continuous on $\mathcal{F}$,
and for some $K>0$, all $f \in\mathcal{F}$,
and all $\varepsilon>0$,
\end{longlist}
\[
\Bigl[\sup_{u>0}u^{2}\Pr^{*}\Bigl(\sup_{\{g\dvtx d_G (g,f)< \varepsilon\}
}|f-g|>u\Bigr)\Bigr]^{1/2}\le K\varepsilon.
\]
Then $\mathcal{F} \in \operatorname{CLT}(P)$.\vadjust{\goodbreak}
\end{theorem}

In this paper we take i.i.d. copies $\{X_{j}\}_{j=1}^{\infty}$ of a
process $\{X(t)\dvtx t \in E \}$, and
consider
%
%
\begin{equation}\label{sets}
\Biggl\{\frac{1}{\sqrt{n}}\sum_{j=1}^n \bigl[I_{X_j(t)\le y}-\Pr
\bigl(X(t)\le y\bigr)\bigr]\dvtx t\in E, y \in\mathbb{R}\Biggr\}
\end{equation}
with the goal of determining when these processes converge in
distribution in some uniform sense to a mean zero Gaussian process. For
example, if the process $X$ has continuous sample paths on $E$, then
$S=C(E)$ and the class of sets $\mathcal{C}$ in (i) of Theorem \ref
{delta-c} consists of the sets $\{f\in C(E)\dvtx f(t)\le y\}$ for $t
\in E$
and $y\in\mathbb{R}$, and we examine when $\mathcal{C} \in \operatorname{CLT}(P)$. As
a result of the previous definitions, the limiting centered Gaussian process
has a version with sample paths in a separable subspace of $\ell
_{\infty
}(\mathcal{C})$, and as a consequence of the addendum to Theorem 1.5.7
of van~der Vaart and Wellner [(\citeyear{vw}), page 37], almost all sample paths are uniformly $L_2$
continuous on $\mathcal{C}$ provided we identify the indicator
functions of sets in $\mathcal{C}$ with $\mathcal{C}$ itself. Furthermore,
we also have that $\mathcal{C}$ is totally bounded in the $L_2$
distance with this identification.
In addition, the following remark is important in this setting.
%
%
\begin{rem}\label{rem1}
The assumption that a centered process $\{X(t)\dvtx t \in T\}$ with $L_2$
distance $d$ is sample bounded and uniformly continuous on $(T,d)$ is
easily seen to follow if $(T,d)$ is
totally bounded and the process is uniformly continuous on $(T,d)$.
Moreover, if the process is Gaussian, the converse also holds using
Sudakov's result as presented in Corollary 3.19 of \citet{led-tal-book}.
\end{rem}

To state and prove our result, we will make use of a distributional
transform that appears in a number of places in the literature; see
\citet{ferguson-book}. \citet{ruschendorf} provides an excellent
introduction to its history, and some uses. In particular, it is used
there to obtain an elegant proof of Sklar's theorem [see \citet{sklar}],
and also in some related applications.

Given the distribution function $F$ of a real valued random variable
$Y$, let $V$ be a random variable uniformly distributed on $[0,1]$ and
independent of~$Y$. In this paper we use the distributional transform
of $F$ defined as
\[
\tilde{F}(x,V)=F(x^{-})+V\bigl(F(x)-F(x^{-})\bigr),
\]
and Proposition 1 in \citet{ruschendorf} shows that
%
%
\begin{equation}\label{uniform}
\tilde{F}(Y,V) \mbox{ is uniform on $[0,1]$}.
\end{equation}
R\"uschendorf calls $\tilde{F}(Y,V)$ the distributional transform of
$Y$, and we also note that $\tilde F(x,V)$ is nondecreasing in $x$.

\section{The main conditions}\label{sec3}

Let $\{X(t)\dvtx t \in E\}$ be a stochastic process as in (\ref{sets}),
and assume
\[
\rho(s,t)=\bigl(\Ex(H_t-H_s)^2\bigr)^{{1}/{2}},\qquad s,t \in E,\vadjust{\goodbreak}
\]
where $\{H(t)\dvtx t \in E\}$ is some Gaussian process which is sample
bounded and $\rho$ uniformly continuous on $E$. In our main result (see
Theorem~\ref{mainresult} below), we hypothesize the relationship between $\{X(t)\dvtx t
\in E\}$ and $\rho(s,t), s,t \in E$ given in (\ref{Lcondition}). The
importance of this condition in the proof of our theorem is 2-fold.
First, it allows one to establish the limiting Gaussian process for our
CLT actually exists. This verifies condition (ii) in Theorem~\ref{agoz} above,
and is accomplished via a subtle application of the necessary and
sufficient conditions for the existence and sample function regularity
of a Gaussian process given in \citet{talagrand-generic}. Second, it
also allows us to verify that the remaining nontrivial condition
sufficient for our CLT, namely condition (iii) of Theorem~\ref{agoz}, applies in
this setting.
This is useful in applications as condition (iii) is in terms of a
single random element involved in our CLT, and hence is far easier to
verify than the typical condition which depends on the full sample of
the random elements as in (\ref{delta-condition}).
The reader should also note that Theorem~\ref{agoz} is a refinement of a number
of previous results sufficient for the CLT, and covers a number of
important earlier empirical process papers.

The existence of such a $\rho(\cdot,\cdot)$
is obtained for a number of specific processes
$\{X(t)\dvtx t \in E\}$ in Section~\ref{sec7}. Nevertheless, given the
process $X$,
determining whether a suitable $\rho$ satisfying (\ref{Lcondition})
exists, or does not exist, may be quite difficult. However, using (\ref
{L1distance})
we may limit our choice for $\rho$ in (\ref{Lcondition}) to be such that
\[
\rho(s,t) \ge c^{-1} \sup_{x \in\mathbb{R}} \bigl(\Ex( |I_{X_{t}\le x}
-I_{X_{s}\le x}|^2)\bigr)^{{1}/{2}}
\]
for all $s,t \in E$ and some constant $c \in(0,\infty)$.

Throughout, $\rho(s,t), s,t \in E$, denotes the $L_2$ metric of a
Gaussian process indexed by $E$ and, to simplify notation, we let
\[
\tilde{F}_{t}(x) \equiv\tilde{(F_t)}(x,V),\qquad x \in\mathbb{R},
\]
be the distributional transform of
$F_t$, the distribution function of $X_t$. Note that this
simplification of notation also includes using $X_t$ for $X(t)$ when
the extra parenthesis make the latter clumsy or unnecessary. Moreover,
this variable notation is also employed for other stochastic processes
in the paper.
In addition, for each
$\varepsilon>0$, let
%
%
\begin{eqnarray}\label{weakLcondition}
&&\sup_{\{s,t \in E\dvtx\rho(s,t)\le
\varepsilon\}}{\Pr}^{*}\bigl(|\tilde{F}_{t}(X_{s}) -\tilde
{F}_{t}(X_{t})|>\varepsilon^{2}\bigr)\le L\varepsilon^{2}\\
&&\qquad(\mbox{the weak $L$ condition})\nonumber
\end{eqnarray}
and
%
%
\begin{eqnarray}\label{Lcondition}
&&\sup_{t \in E}\Pr^{*}\Bigl({\sup_{\{
s\dvtx\rho
(s,t)\le\varepsilon\}}}|\tilde{F}_{t}(X_{s}) -\tilde
{F}_{t}(X_{t})|>\varepsilon^{2}\Bigr)\le L\varepsilon^{2}\\
&&\qquad(\mbox{the $L$ condition}).\nonumber
\end{eqnarray}

%
\begin{rem}\label{Cepsilon^{2}}In the $L$ conditions the probabilities
involve an $\varepsilon^{2}$. However, since for any constant $C\in
(0,\infty)$, an $L_{2}$ metric $\rho$, is pre-Gaussian if and only if
$C\rho$ is pre-Gaussian, WLOG we can change to $C\varepsilon^{2}$.
Moreover, note that any constant $L$ sufficient for (\ref{Lcondition})
will also suffice for (\ref{weakLcondition}), and hence to simplify
notation we do not distinguish between them.
\end{rem}
%
%
\begin{lem} \label{weakLlemma} Let $L$ be as in (\ref{weakLcondition}),
and take $s,t \in E$. Then, for all $ x \in\mathbb{R}$,
%
%
\begin{equation}\label{basicLineq}
\Pr(X_s \leq x < X_t) \leq(L+1)\rho^2(s,t)
\end{equation}
and by symmetry,
%
%
\begin{eqnarray} \label{L1distance}
\Ex|I_{X_{t}\le x} -I_{X_{s}\le x}|&=&\Pr(X_t \leq x <
X_s)+\Pr(X_s \leq x < X_t)\nonumber\\[-8pt]\\[-8pt]
&\le&2(L+1)\rho^2(s,t).\nonumber
\end{eqnarray}
Further, we have
%
%
\begin{equation}\label{weakLdistribution}
{\sup_{x}}|F_{t}(x)-F_{s}(x)|\le
2(L+1)\rho^2(s,t).
\end{equation}
\end{lem}
%
%
\begin{rem} As in Lemma~\ref{weakLlemma} and Lemmas~\ref{tauestimate} and
\ref{taulower} below, use only the weak $L$
condition (\ref{weakLcondition}). Actually, for Lemmas
\ref{tauestimate} and~\ref{taulower}, all we
need is Lemma~\ref{weakLlemma}. However, in Lemma~\ref{lem4} we need the stronger form as
stated in (\ref{Lcondition}).
\end{rem}
\begin{pf*}{Proof of Lemma~\ref{weakLlemma}}
Since $\tilde F_t$ is nondecreasing and $x<y$ implies
$F_t(x) \leq\tilde F_t(y)$, we have
\[
\Pr(X_s \leq x < X_t) \leq\Pr\bigl(\tilde F_t(X_s) \leq\tilde
F_t(x),F_t(x) \leq\tilde F_t(X_t)\bigr).
\]
Thus
\begin{eqnarray*}
\Pr(X_s \leq x < X_t) &\leq&\Pr\bigl(F_t(x) \leq\tilde F_t(X_t) \leq F_t(x)
+ \rho^2(s,t),\tilde F_t(X_s) \leq\tilde F_t(x)\bigr)\\
&&{} +
\Pr\bigl(\tilde F_t(X_t) > F_t(x) + \rho^2(s,t),\tilde F_t(X_s) \leq
\tilde F_t(x)\bigr)
\end{eqnarray*}
and hence
\begin{eqnarray*}
\Pr(X_s \leq x < X_t) &\leq&\Pr\bigl(F_t(x) \leq\tilde F_t(X_t) \leq F_t(x)
+ \rho^2(s,t)\bigr)
\\
&&{} + \Pr\bigl(|\tilde F_t(X_t) - \tilde F_t(X_s)|
>\rho
^2(s,t)\bigr).
\end{eqnarray*}
Now (\ref{weakLcondition}) implies for all $s,t \in E$ that
\[
\Pr\bigl(|\tilde F_t(X_t) - \tilde F_t(X_s)| >\rho^2(s,t)\bigr)\le L\rho^2(s,t),
\]
since its failure for $s_0,t_0 \in E$ and $\varepsilon= \rho(s_0,t_0)$ in
(\ref{weakLcondition}) implies a contradiction. Therefore, since
$\tilde F_t (X_t)$ is uniform on $[0,1]$, we 
have
\[
\Pr(X_s \leq x < X_t) \leq\rho^2(s,t) + L\rho^2(s,t).
\]
The last conclusion follows by moving the absolute values outside the
expectation.
\end{pf*}

\section{The main result}\label{sec4}

Recall the relationship of the $X$ process and $\rho$ as described at
the beginning of Section~\ref{sec3}. Then we have:
%
%
\begin{theorem}\label{mainresult} Let $\rho$ be given by $\rho
^{2}(s,t)=\Ex
(H(s)-H(t))^{2}$, for some centered Gaussian process $H$ that is sample
bounded and uniformly continuous on $(E,\rho)$ with probability one.
Furthermore, assume that for some $L<\infty$, and all $\varepsilon>0$,
the $L$ condition (\ref{Lcondition}) holds, and $D(E)$ is a collection
of real valued functions on $E$ such that $\Pr(X(\cdot) \in D(E))=1$.
If
\[
\mathcal{C}=\{C_{s,x}\dvtx s \in E, x\in\mathbb{R}\},
\]
where
\[
C_{s,x}=\{z\in D(E)\dvtx z(s) \le x\}
\]
for $s \in E,x \in\mathbb{R}$, then $\mathcal{C} \in \operatorname{CLT}(P)$.
\end{theorem}
%
%
\begin{rem}
Note that a sample function of the $X$-process is in the set $C_{s,x}$
iff $X_s \le x$. Hence, if we
identify a point $(s,x) \in E \times\mathbb{R}$ with the set
$C_{s,x}$, then instead of saying $\mathcal{C} \in \operatorname{CLT}(P)$, we will
often say
\[
\{I_{X_{s}\le x}-\Pr(X_{s}\le x)\dvtx s\in E, x\in\mathbb{R}\}
\]
satisfies the CLT in $\ell_{\infty}(\mathcal{C})$ [or in $\ell
_{\infty
}(E \times\mathbb{R})$].
\end{rem}
%
%
\begin{rem}
At this point one might guess that the reader is questioning the
various assumptions in Theorem~\ref{mainresult}. First we mention that
$D(E)$ is some convenient function space. For example, typically the
process $X$ has continuous sample paths on $E$, so $D(E)=C(E)$ in these
situations. More perplexing, at least for most readers, is probably the
appearance of the distributional transforms $\{\tilde{F}_{t}\dvtx t \in
E\}$ in the $L$ condition (\ref{Lcondition}). If the distribution
functions $F_t$ are all continuous, then $F_t= \tilde{F}_{t}, t \in E$,
and our proof obviously holds with $F_t$ replacing $\tilde{F}_{t}$ in
the $L$ condition. However, without all the distribution functions
$F_t$ assumed continuous, the methods required in our proof fail with
this substitution. An interesting case where the distributional
transforms are useful occurs when one has a point $t_0 \in E$ such that
$\Pr (X(t)=X(t_0) \mbox{ for all } t\mbox{ in }E)=1$, and $F_{t_0}$ is
possibly discontinuous. In this situation, the $ L$ condition
(\ref{Lcondition}) holds for the Gaussian process $H(t)=g$ for all $ t
\in E$, $ g$ a standard Gaussian random variable and $X(t_0)$ having
any distribution function $F_{t_0}$. Thus Theorem~\ref{mainresult}
applies and yields the classical empirical CLT when the set $S$ is the
real line, and the class of sets consists of half-lines for all laws
$F_{t_0}$. A similar result also applies if $E$ is a finite disjoint
union of nonempty sets, and the process $\{X_t\dvtx t \in E\}$ is
constant on each of the disjoint pieces of $E$ regardless of the
distribution functions $F_t, t \in E$. More importantly, however,
allowing even a single discontinuous distribution $F_t$ may invalidate
the empirical $\mathrm{CLT}$ on $\mathcal{C}$. For example, if
$\{X(t)\dvtx t \in[0,1]\}$ is standard Brownian motion with
$P(X(0)=0)=1$, then in Section \ref {sec81} we show the empirical
$\mathrm{CLT}$ fails, but Corollary~\ref{cor2} shows that it holds if
we allow the distribution at time zero to have a bounded density.
Furthermore, in Section~\ref{sec84} we provide some additional examples
where the empirical process is pre-Gaussian, and the input process
$\{X(t)\dvtx t \in E\}$ satisfies a modified $L$~condition, that is,
for all $\varepsilon>0$, there is an $L<\infty$ such that
\[
\sup_{t \in E} \Pr^{*}\Bigl({\sup_{\{s\dvtx\rho(s,t) \leq\varepsilon\}
}}|F_t(X_s)-F_t(X_t)|> \varepsilon^2\Bigr) \leq\varepsilon^2,
\]
yet the empirical $\mathrm{CLT}$ we seek fails. Hence one needs to assume
something more, and our results show that the $L$ condition given in
(10) is sufficient for the empirical~$\mathrm{CLT}$.
\end{rem}

In Section~\ref{sec7} we will provide another theorem showing how
Theorem~\ref{mainresult} can be applied, and hence the examples
obtained there are motivation for its formulation in terms of the $L$
condition (\ref{Lcondition}). The following remark also motivates the
presence of the process $\{\tilde{F}_{t}(X_{s})\dvtx s \in E\}$ and the
$L$ condition in our CLT. In particular, we sketch an argument that for
each $t \in E$ a symmetric version of this process satisfies a CLT in
$\ell_{\infty}(E)$. This remark is meant only for motivation, and in
its presentation we are unconcerned with a number of details.
%
%
\begin{rem}\label{NecCond}Let
\[
\{I_{X_{s}\le x}-\Pr(X_{s}\le x)\dvtx s\in E, x\in\mathbb{R}\}
\]
satisfy the central limit theorem in the closed subspace of $\ell
_{\infty}(E \times\mathbb{R})$ consisting of functions whose
s-sections are Borel measurable on $\mathbb{R}$. We denote this
subspace by
$\ell_{\infty, m}(E\times\mathbb{R})$, and also assume the
distribution functions $F_t$ are all continuous.
Then, for each fixed $t \in E$, we define the bounded linear operator
$\phi\dvtx\ell_{\infty,m}(E\times\mathbb{R})\lra\ell_{\infty}(E)$
given by
\[
\phi(f)(s)=\int f(s,x)F_{t}(dx).
\]
Now by the symmetrization lemma [Lemma 2.7 in \citet
{empirical-specialinvited}], we have for a Rademacher random variable
$\varepsilon$ independent of the empirical process variables that
$\{\varepsilon I_{C_{s,x}}\dvtx s \in E, x \in\mathbb{R}\}$ satisfy the CLT
in $\ell_{\infty}(E \times\mathbb{R})$.
Taking
$f(s,x)= I_{C_{s,x}}$, we have for all $t \in E$ fixed that
\[
\phi(f)(s)= 1 - F_t(z(s)^{-})= 1 - F_t(z(s))
\]
as we are assuming the $F_t$ are continuous. Therefore the continuous
mapping theorem [see, e.g., Theorem 1.3.6 in \citet{vw}] implies that
for each $t\in E$,
\[
Z_n(s)=\frac1{\sqrt{n}}\sum_{j=1}^{n}\varepsilon_{j}\bigl(1- F_t(X_s)\bigr),
\qquad s\in E,\vadjust{\goodbreak}
\]
satisfies the CLT in $\ell_{\infty}(E)$. In addition,
since we are assuming $F_t=\tilde F_t$, we should then have
``asymptotically small oscillations;'' namely, for every $\delta>0$
there exists $\varepsilon>0$ such that
\[
{\Pr}^*\Biggl(\sup_{\rho(s,t)\le\varepsilon}\frac1{\sqrt
{n}}\sum_{j=1}^{n}\varepsilon_{j}\bigl(\tilde{F}_{t}(X_{s})-\tilde
{F}_{t}(X_{t})\bigr)>\delta\Biggr)\le\delta.
\]
By using standard symmetry arguments this last probability dominates
\[
\frac12 \Pr^* \Bigl({\max_{j\le n}\sup_{\rho(s,t)\le
\varepsilon}}|\tilde{F}_{t}(X_{s})-\tilde{F}_{t}(X_{t})|>\sqrt{n}\delta
\Bigr)\le\delta,
\]
which (again by standard arguments) implies (modulo multiplicative constants)
\[
n\Pr^* \Bigl({\sup_{\rho(s,t)\le\varepsilon}}|\tilde
{F}_{t}(X_{s})-\tilde{F}_{t}(X_{t})|>\sqrt{n}\delta\Bigr)\le\delta.
\]
While this is different from the hypotheses in our theorem, it
indicates that the quantity ${\sup_{\rho(s,t)\le\varepsilon}}|\tilde
{F}_{t}(X_{s})-\tilde{F}_{t}(X_{t})|$ is relevant to any such theorem.
\end{rem}
%

\section{Preliminaries for generic chaining}\label{sec5}

Let $T$ be an arbitrary countable set. Then, following \citet
{talagrand-generic} we have:
%
%
\begin{mydef}
An admissible sequence is an increasing sequence
$(\mathcal
{A}_{n})$ of partitions of $T$ such that
\[
\card\mathcal{A}_{n}\le N_{n},
\]
where $N_0=1$, and for $n \geq1$, $N_n=2^{2^{n}}$.
The partitions $(\mathcal{A}_{n})$ are increasing if every set in
$\mathcal{A}_{n+1}$ is a subset of some set of $\mathcal{A}_n$.
\end{mydef}

We also have:
%
%
\begin{mydef}
If $t \in T$, we denote by $A_n(t)$ the unique element of $\mathcal
{A}_n$ that contains $t$. For a psuedo-metric $e$ on $T$, and $A
\subseteq E$, we write $\Delta_{e}(A)$ to denote the diameter of $A$
with respect to $e$.
\end{mydef}

Using generic chaining and the previous definitions, Theorem 1.4.1 of
\citet{talagrand-generic} is essentially the following result. Its
statement there contains a curious wording at the end of the first
sentence, which suggests that cutting and pasting led to something
being omitted. After closer inspection we observed that the necessary
assumption of total boundedness of the parameter space was required,
and it now appears in the statement of the theorem below. Since Theorem
1.4.1 appears without proof, for completeness the proof can be found in
the \hyperref[app]{Appendix}.
%
%
\begin{theorem} \label{talagrand-generic}Let $\{X_t\dvtx t \in T\}$ be a
centered Gaussian process with $L_2$ distance $d(s,t), s,t \in T$,
where $T$ is countable, and $(T,d)$ is totally bounded. Then, the
following are equivalent:\vadjust{\goodbreak}

\begin{longlist}
\item
$X_t$ is uniformly continuous on $(T, d)$ with probability one.

\item We have
\[
\lim_{\varepsilon\rightarrow0} \Ex\Bigl(\sup_{d(s,t) \le\varepsilon}(X_s -
X_t)\Bigr) = 0.
\]

\item There exists an admissible sequence of partitions of $T$ such that
%
%
\begin{equation}\label{Talagrand}
\lim_{k \rightarrow\infty} \sup_{t \in T} \sum_{n \geq k}
2^{n/2}\Delta
(A_n(t)) =0.
\end{equation}
\end{longlist}
\end{theorem}

Under the assumption that $H$ is centered Gaussian and uniformly
continuous on $(T,e)$, then, recalling Remark~\ref{rem1}, it follows that $H$
being sample bounded on $T$ is equivalent to $(T,e)$ being totally
bounded. Also, an immediate corollary of this result used below is as follows.
%
%
\begin{prop}\label{prop1} Let $H_1$ and $H_2$ be mean zero Gaussian processes with
$L_2$ distances $e_1, e_2$, respectively, on $T$. Furthermore, assume T
is countable, and $e_1(s,t) \leq e_2(s,t)$ for all $s,t \in T$.
Then, $H_2$ sample bounded and uniformly continuous on $(T,e_2)$ with
probability one, implies $H_1$ is sample bounded and uniformly
continuous on $(T,e_1)$ with probability one.
\end{prop}
%
%
\begin{rem} One can prove this using Slepian's lemma [see, e.g., \citet
{fernique-1975}]. However, the immediate conclusion is that $H_{1}$ is
sample bounded and uniformly continuous on $(T,e_{2})$. Then, a
separate argument is needed to show the statement in this proposition.
Using the more classical formulation for continuity of Gaussian
processes involving majorizing measures [see, e.g., Theorem~12.9 of
\citet{led-tal-book}], the result also follows similarly to what is
explained below.
\end{rem}
\begin{pf*}{Proof of Proposition~\ref{prop1}}
By the previous theorem $\{H_2(t)\dvtx t \in T\}$ is sample
bounded and uniformly continuous
on $(T,e_2)$ with probability one if and only if there exists an
admissible sequence of partitions of $T$ such that
\[
\lim_{r \rightarrow\infty} \sup_{t \in T} \sum_{n \geq r} 2^{n/2}
\Delta_{e_2}(A_n(t))=0.
\]
Since $\Delta_{e_1}(A_n(t)) \leq\Delta_{e_2}(A_n(t))$, we have
\[
\lim_{r \rightarrow\infty} \sup_{t \in T} \sum_{n \geq r} 2^{n/2}
\Delta_{e_1}(A_n(t))=0,
\]
and hence Theorem~\ref{talagrand-generic} implies that $H_1$ is sample bounded and uniformly
continuous on $(T,e_1)$ with probability one. Thus the proposition is proven.
\end{pf*}

\section{\texorpdfstring{Proof of Theorem \protect\ref{mainresult}}{Proof of Theorem 3}}\label{sec6}

First we establish some necessary lemmas, and the section ends with the
proof of Theorem~\ref{mainresult}. Throughout we take as given the assumptions and
notation of that theorem.\vadjust{\goodbreak}

\subsection{Some additional lemmas}\label{sec61}

In order to simplify notation, we denote the $L_2$ distance on the
class of indicator functions
\[
\mathcal{F}= \{I_{X_s \leq x}\dvtx s \in E, x \in\mathbb{R}\}
\]
by writing
\[
\tau((s,x),(t,y)) = \bigl\{E\bigl(( I_{X_s \le x} - I_{X_t \le y})^2\bigr)\bigr\}^{{1}/{2}}
\]
and identifying $\mathcal{F}$ with $E \times\mathbb{R}$.
Our next lemma relates the $\tau$-distance and the $\rho$-distance. It
upgrades (\ref{L1distance}) when $x \not=y$.
%
%
\begin{lem}\label{tauestimate}Assume that (\ref{weakLcondition})
holds. Then
%
%
\begin{equation}\label{tau-rho}
\tau^{2}((s,x),(t,y)) \leq{\min_{u \in\{s,t\}}}|F_{u}(y)-F_{u}(x)| +
(2L+2)\rho^{2}(t,s).
\end{equation}
Moreover, if $Q$ denotes the rational numbers, there is a countable
dense set $E_0$ of $(E,\rho)$ such that $\mathcal{F}_0= \{I_{X_s \le
x}\dvtx(s,x) \in E_0 \times Q\}$ is dense in
$(\mathcal{F},\tau)$.
\end{lem}
\begin{pf}
First observe that by using the symmetry in $s$ and $t$
of the right-hand term of (\ref{tau-rho}), we have, by applying (\ref
{basicLineq}) in the second inequality below, that
\begin{eqnarray*}
\tau^{2}((s,x),(t,y))&=&\Ex|I_{X_{t}\le y}-I_{X_{s}\le x}|\le\Ex
|I_{X_{t}\le y}-I_{X_{t}\le x}|+\Ex|I_{X_{t}\le x}-I_{X_{s}\le x}|
\\
&=&|F_{t}(y)-F_{t}(x)| + \Pr(X_{s}\le x<X_{t}) + \Pr
(X_{t}\le x<X_{s})
\\
&\leq& |F_{t}(y)-F_{t}(x)| + (2+2L)\rho^2(s,t).
\end{eqnarray*}
Similarly, we also have by applying (\ref{basicLineq}) again that
\[
\tau^{2}((s,x),(t,y))
\leq|F_{s}(y)-F_{s}(x)| + (2+2L)\rho^2(s,t).
\]
Combining these two inequalities for $\tau$, the proof of (\ref
{tau-rho}) holds.
Since $(E, \rho)$ is assumed totally bounded, there is a countable
dense set $E_0$ of $(E,\rho)$, and hence the right continuity of the
distribution functions and (\ref{tau-rho}) then imply the final
statement in Lemma~\ref{tauestimate}.
\end{pf}

Using Lemma~\ref{tauestimate} and the triangle inequality, we can estimate the $\tau
$-diameter of sets as follows.
%
%
\begin{cor}\label{taudiam}If $t_{B}\in B\subseteq E$ and $D\subseteq
\mathbb{R}$, then
\[
\diam_{\tau}(B\times D)\le2\Bigl\{ (2L+2)^{1/2}\diam_{\rho}(B) + {\sup
_{x,y\in D}}|F_{t_{B}}(y)-F_{t_{B}}(x)|^{1/2}\Bigr\}.
\]
\end{cor}
%
%
\begin{lem}\label{taulower} Assume that $(s,x)$ and $(t,y)$ satisfy
\[
\tau((s,x),(t,y))=\|I_{X_{s}\le x}-I_{X_{t}\le y}\|_{2}\le\varepsilon,
\]
$\rho(s,t)\le\varepsilon$, and (\ref{weakLcondition}) holds. Then, for
$c=(2L+2)^{1/2}+1$,
\[
|F_{t}(x) - F_{t}(y)| \le(c\varepsilon)^{2}\vadjust{\goodbreak}
\]
or, in other words,
\[
|F_{t}(x) - F_{t}(y)| \le(c\max\{\tau((s,x),(t,y)),\rho(s,t)\})^{2}.
\]
\end{lem}
\begin{pf} Using (\ref{basicLineq}) in the second inequality below,
we have
\begin{eqnarray*}
|F_{t}(y)-F_{t}(x)|^{1/2}
&=&\|I_{X_{t}\le x}-I_{X_{t}\le y}\|_{2}\le\|I_{X_{s}\le x}-I_{X_{t}\le
y}\|_{2}+\|I_{X_{s}\le x}-I_{X_{t}\le x}\|_{2}
\\
&\le&\varepsilon+\bigl(\Pr(X_{s}\le x<X_{t})+\Pr
(X_{t}\le x< X_{s})\bigr)^{1/2}
\\
&\leq&\varepsilon+(2L\varepsilon^{2}+2\varepsilon
^{2})^{1/2}=[(2L+2)^{1/2}+1]\varepsilon\equiv c\varepsilon.
\end{eqnarray*}
Hence the lemma is proven.
\end{pf}

The next lemma is an important step in verifying the weak-$L_2$
condition in item (iii) of Theorem~\ref{agoz} above; for example, see Theorem
4.4 of \citet{clt-localcondit}.
%
%
\begin{lem}\label{lem4} If (\ref{Lcondition}) holds, $c$ is as in Lemma~\ref{taulower}, and
\[
\lambda((s,x),(t,y)) = \max\{ \tau((s,x),(t,y)), \rho(s,t)\};
\]
then for all $(t,y)$ and $\varepsilon>0$,
\[
\Pr^{*}\Bigl({\sup_{\{(s,x)\dvtx\lambda((t,y),(s,x))\le\varepsilon\}
}}|I_{X_{t}\le y} - I_{X_{s}\le x}|>0\Bigr)\le2(c^{2}+L+1)\varepsilon^{2}.
\]
\end{lem}
\begin{pf}
First we observe that
\begin{eqnarray*}
&&{\Pr}^*\Bigl({\sup_{\{(s,x)\dvtx\lambda((t,y),(s,x))\le
\varepsilon\}
}}|I_{X_{t}\le y} -I_{X_{s}\le x}|>0\Bigr)\\
&&\qquad= {\Pr}^*\Bigl(\sup_{\{(s,x)\dvtx\lambda((t,y),(s,x))\le\varepsilon\}
}I_{X_{t}\le
y,X_{s}>x} + I_{X_{s}\le x,X_{t}>y}>0\Bigr).
\end{eqnarray*}
Again, using the fact that $x<y$ implies $F_t(x) \leq\tilde F_t(y)$,
we have
\begin{eqnarray*}
&&{\Pr}^*\Bigl({\sup_{\{(s,x)\dvtx\lambda((t,y),(s,x))\le
\varepsilon\}
}}|I_{X_{t}\le y} -I_{X_{s}\le x}|>0\Bigr)\\
&&\qquad\le{\Pr}^*\Bigl(\sup_{\{(s,x)\dvtx\lambda((t,y),(s,x))\le\varepsilon\}
}I_{\tilde
{F}_{t}(X_{t})\le\tilde{F}_{t}(y), {F}_{t}(x)\le\tilde
{F}_{t}(X_{s})}>0\Bigr)\\
&&\qquad\quad{}+{\Pr}^*\Bigl(\sup_{\{(s,x)\dvtx\lambda((t,y),(s,x))\le
\varepsilon\}}I_{\tilde{F}_{t}(X_{s})\le\tilde{F}_{t}(x),
{F}_{t}(y)\le
\tilde{F}_{t}(X_{t})}>0\Bigr)=I + \mathit{II},
\end{eqnarray*}
where
%
%
\begin{equation}\label{inf}
I= {\Pr}^*\Bigl(\sup_{\{(s,x)\dvtx\lambda((t,y),(s,x))\le\varepsilon\}
}I_{\tilde
{F}_{t}(X_{t})\le\tilde{F}_{t}(y), {F}_{t}(x)\le\tilde
{F}_{t}(X_{s})}>0\Bigr)
\end{equation}
and
%
%
\begin{equation}\label{I,II}
\mathit{II}={\Pr}^*\Bigl(\sup_{\{(s,x)\dvtx\lambda((t,y),(s,x))\le\varepsilon\}
}I_{\tilde
{F}_{t}(X_{s})\le\tilde{F}_{t}(x), {F}_{t}(y)\le\tilde
{F}_{t}(X_{t})}>0\Bigr).\vadjust{\goodbreak}
\end{equation}
At this point we use Lemma~\ref{taulower} to see that
in (\ref{inf}) we can use
\[
\inf_{\{(s,x)\dvtx\lambda((t,y),(s,x))\le\varepsilon\}}F_{t}(x)\ge
F_{t}(y)-(c\varepsilon)^{2}.
\]
Therefore, since $\tilde F_t(x) \leq F_t(x)$ for all $x$ and again
using (\ref{Lcondition})
\begin{eqnarray*}
I&\le&{\Pr}^*\Bigl(\sup_{\{(s,x)\dvtx\lambda((t,y),(s,x))\le
\varepsilon\}}I_{\tilde{F}_{t}(X_{t})\le{F}_{t}(y),
F_{t}(y)-(c\varepsilon
)^{2}\le\tilde{F}_{t}(X_{s})}>0\Bigr)\\
&\le&{\Pr}^*\Bigl(\sup_{\{(s,x)\dvtx\lambda((t,y),(s,x))\le\varepsilon\}
}I_{\tilde
{F}_{t}(X_{t})\le{F}_{t}(y), F_{t}(y)-(c\varepsilon)^{2}\le\tilde
{F}_{t}(X_{t})+ \varepsilon^2 }
>0\Bigr) +L\varepsilon^{2}\\
&\le&\Pr\bigl({F_{t}(y)-(c\varepsilon)^{2}-\varepsilon^{2}\le\tilde
{F}_{t}(X_{t})\le{F}_{t}(y)}\bigr)+L\varepsilon^{2}\\
&\le&(c^{2}+L+1)\varepsilon^{2} \qquad\mbox{by (\ref{uniform})}.
\end{eqnarray*}
Now, we estimate II in (\ref{I,II}). Again using the fact that $\tilde
F_t(x) \leq F_t(x)$ for all $x$, Lemma~\ref{taulower}, and our definition of $L$, we
therefore have
\begin{eqnarray*}
&&{\Pr}^{*}\Bigl(\sup_{\{(s,x)\dvtx\lambda((t,y),(s,x))\le
\varepsilon
\}}I_{\tilde{F}_{t}(X_{s})\le\tilde{F}_{t}(x), {F}_{t}(y)\le\tilde
{F}_{t}(X_{t})}>0\Bigr)\\
&&\qquad\le\Pr\bigl(\tilde{F}_{t}(X_{t})-\varepsilon^{2}\le
F_{t}(y)+(c\varepsilon
)^{2}, {F}_{t}(y)\le\tilde{F}_{t}(X_{t})\bigr)+L\varepsilon^{2}\\
&&\qquad\le(c^{2}+L+1)\varepsilon^{2}.
\end{eqnarray*}
\upqed
\end{pf}

\subsection{\texorpdfstring{The construction and the proof of Theorem \protect\ref{mainresult}}{The construction and the proof of Theorem 3}}\label{sec62}
Since $(E, \rho)$ is totally bounded by Remark~\ref{rem1}, take $E_0$ to be any
countable dense subset of $E$ in the $\rho$ distance. Then by Theorem
\ref{talagrand-generic}, Talagrand's continuity theorem, there exists
an admissible sequence of partitions, $\mathcal{B}_{n}$ of $E_0$, for which
%
%
\begin{equation}\label{gamma-tau}
\lim_{r\rightarrow\infty}\sup_{t \in
E_0}\sum_{n\ge r}2^{n/2}\Delta_{\rho}(B_{n}(t)) =0.
\end{equation}

Fix $n$. Then, for each
$B\in\mathcal{B}_{n-1}$ choose $t_{B}\in B$.
Fix the distribution function $F_{B}:=F_{t_{B}} $ and $\mu_{B}$ the
associated probability measure.
Put $\alpha=(\Delta_{\rho}(B) + 2^{-n})^{2}$ and set $z_{1}=\sup\{
x\in
\mathbb{R}\dvtx F_{B}(x) < \alpha\}$. We consider two cases:
\begin{itemize}
\item$F_{B}(z_{1})\le\alpha$ and
\item$F_{B}(z_{1})> \alpha$.
\end{itemize}

In the first case $F_{B}(z_{1})=\alpha$. If $F_{B}(z_{1})<\alpha$, then
by right continuity there exist $w>z_{1}$ such that $F_{B}(w)<\alpha$,
which contradicts the definition of $z_{1}$.
In this case we consider $C_{1}=(-\infty, z_{1}]$ and
$D_{1}=\varnothing$.

In the second case we let $C_{1}=(-\infty, z_{1})$ and $D_{1}=\{z_{1}\}
$. In either case $\mu_B(C_{1}\cup D_{1})\ge\alpha$.

If $\mu_{B}((z_{1},\infty))\ge\alpha$, let $z_{2}=\sup
\{
x>z_{k}\dvtx F_{B}(x)-F_{B}(z_{1})< \alpha\}$.
If $z_{2}=\infty$, we set $C_{2}=(z_1,\infty)$ and $D_{2}=
\varnothing$.
Otherwise, if $z_{2} < \infty$, there are two cases.
That is, we have:
\begin{itemize}
\item$F_{B}(z_{2})-F_{B}(z_{1})\le\alpha$ and
\item$F_{B}(z_{2})-F_{B}(z_{1})> \alpha$.\vadjust{\goodbreak}
\end{itemize}
In the first case we consider $C_{2}=(z_{1}, z_{2}]$ and
$D_{2}=\varnothing$. In the second case we let $C_{2}=(z_{1}, z_{2})$ and
$D_{2}=\{z_{2}\}$. As before, $\mu_B(C_{2}\cup D_{2})\ge\alpha$.

Now assume that we have constructed $C_{1},\ldots,C_{k}$ and
$D_{1},\ldots,D_{k}$ in this manner. Therefore we have $z_{k}$. If
$\mu
_{B}((z_{k},\infty))\ge\alpha$, let $z_{k+1}=\sup\{x>z_{k}\dvtx\break
F_{B}(x)-F_{B}(z_{k})< \alpha\}$.
If $z_{k+1}=\infty$, we set $C_{k+1}=(z_k,\infty)$ and $D_{k+1}=
\varnothing$. Otherwise, if $z_{k+1} < \infty$, there are two cases. That
is, we have:
\begin{itemize}
\item$F_{B}(z_{k+1})-F_{B}(z_{k})\le\alpha$ and
\item$F_{B}(z_{k+1})-F_{B}(z_{k})> \alpha$.
\end{itemize}
In the first case we consider $C_{k+1}=(z_{k}, z_{k+1}]$ and
$D_{k+1}=\varnothing$. In the second case we let $C_{k+1}=(z_{k},
z_{k+1})$ and $D_{k+1}=\{z_{k+1}\}$. As before, $\mu_B(C_{k+1}\cup
D_{k+1})\ge\alpha$.
Hence, there can be at most $\frac1{\alpha}+1$ steps before $\{
C_{k},D_{k}\}_{k}$ cover $\mathbb{R}$. Therefore, after eliminating any
empty set, we have a cover of $\mathbb{R}$ with at most
$\frac2{\alpha
}+2$ sets. By our choice of $\alpha$ the cover has at most
$2^{2n+1}+2$ sets.
Hence since we have $B \in\mathcal{B}_{n-1}$, the number of sets
used to cover $E_0 \times\mathbb{R}$ of the form $B \times C_k$ or $B
\times D_k$ is less than or equal to $ 2^{2^{n-1}}(2^{2n+1}+2)$.
The reader should note that the points $\{z_k\}$ depend on the set $B$,
but we have suppressed that to simplify notation.
We now check the $\tau$-diameters of the nonempty $B\times C_{k}$ and
$B\times D_{k}$.

Estimating these diameters by doubling the radius of the sets, the
triangle inequality allows us to upper bound their radius using one of
$s$ and $t$ to be $t_{B}$. Also note that in Lemma~\ref{tauestimate}, or Corollary~\ref{taudiam},
the term which contains $|F_{t_{B}}(y)-F_{t_{B}}(x)|$ would cause
trouble in the case $D_{k}\neq\varnothing$, since this is only known to
be $\ge\alpha$. Luckily it does not appear when $D_k \not=
\varnothing$.

First we consider the $\tau$-diameter of sets of the form $B\times
C_{k}$ when
$D_{k}=\varnothing$. Then $C_k=(z_{k-1,B}, z_{k,B}]$. Hence for
$(s,x), (t,y)\in B\times C_{k}$, Corollary~\ref{taudiam} implies
\[
\Delta_{\tau}\bigl(B\times(z_{B,k-1},z_{B,k}]\bigr)\le2
\biggl((2L+2)^{1/2}\Delta
_{\rho}(B) + \Delta_{\rho}(B) + \frac1{2^{n}}\biggr).
\]
When $D_{k}\neq\varnothing$, then $C_k=(z_{k-1,B}, z_{k,B})$, so
again by
Corollary~\ref{taudiam} the $\tau$-diameter of $B \times C_k$ has an upper bound as
in the previous case.

If $D_{k}\neq\varnothing$, then the only element of $D_{k}$ is $z_{k,B}$,
and by Corollary~\ref{taudiam} we have
\[
\Delta_{\tau}(B \times D_k) \leq2(2L+2)^{{1}/{2}}\diam_{\rho}(B).
\]
So, in either case,
%
%
\begin{equation}\label{Delta}
\Delta_{\tau}(B\times C_{B,k}\mbox{ or
}D_{B,k})\le2\biggl((2L+2)^{1/2}\Delta_{\rho}(B) + \Delta_{\rho
}(B) +
\frac1{2^{n}}\biggr).
\end{equation}

%
\begin{lem}\label{partitions}Let $\mathcal{G}_{n}$ be a sequence of
partitions of an arbitrary parameter set $T$ with pseudo metric $e$ on
$T$ satisfying both:
\begin{longlist}[(2)]
\item[(1)] $\card(\mathcal{G}_{n})\le2^{2^{n}}$ and
\item[(2)] $\lim_{r\to\infty}\sup_{t \in T}\sum_{n\ge
r}2^{n/2}\Delta
_{e}(G_{n}(t))=0$,\vadjust{\goodbreak}
\end{longlist}
and set\vspace*{1pt} $\mathcal{H}_{n}:=\mathcal{P}(\bigcup_{1 \le k\le
n-1}\mathcal
{G}_{k})$, where $\mathcal{P}(\mathcal{D})$ denotes the minimal
partition generated by the sets in $\mathcal{D}$. Then the sequence
$\mathcal{H}_n$ (\textup{notice the $n-1$ in the union}) also satisfies
those conditions.
\end{lem}
\begin{pf}
The first condition holds since a simple induction on $n$
implies the minimal partition
\[
\mathcal{H}_n=\mathcal{P}\biggl( \bigcup_{1 \leq k \leq n-1} \mathcal{G}_k\biggr)
\]
has cardinality at most $\prod_{k=1}^{n-1} 2^{2^{k}} \leq2^{2^n}$.
The second condition holds since the partitions are increasing
collections of sets, and hence
$\diam_{e}(H_{n}(t))\le\diam_{e}(G_{n-1}(t))$.
\end{pf}
%
%
\begin{lem}\label{tauPG}Let $E_0$ be a countable dense subset of $(E,
\rho)$. Then there exists an admissible sequence of partitions $\{
\mathcal{A}_n\dvtx n \ge0\}$ of $E_0 \times\mathbb{R}$ such that
%
%
\begin{equation}\label{tauPreG}
\lim_{r\to\infty}\sup_{(t,y) \in E_0\times\mathbb{R}}\sum_{n\ge
r}2^{n/2}\Delta_{\tau}(A_{n}((t,y)))=0.
\end{equation}
\end{lem}
\begin{pf}
We construct the admissible sequence of partitions
$\mathcal{A}_{n}$ as above. More precisely, let $\{\mathcal{B}_n\dvtx n
\ge
0\}$ be an increasing sequence of partitions of $E_0$ such that (\ref
{gamma-tau}) holds, and after the construction above we also have (\ref
{Delta}).
That is, for $k \ge1$ let
\[
\mathcal{G}_k =\{B \times F\dvtx B \in\mathcal{B}_{k-1}, F \in\mathcal
{E_B}\},
\]
where
\[
\mathcal{E}_B = \{C_{j,B},D_{j,B} \mbox{ all sets nonempty}\}
\]
and $C_{j,B},D_{j,B}$ are constructed from $B \in\mathcal{B}_{k-1}$ as
above. Then, for $n \ge4$ set
\[
\mathcal{A}_n =\mathcal{P}\biggl( \bigcup_{3 \leq k \leq n-1} \mathcal{G}_k\biggr),
\]
where $\mathcal{P}(\mathcal{D})$ is the minimal partition generated by
the sets in $\mathcal{D}$, and for $n=1,2,3$ we take $\mathcal{A}_n$ to
be the single set $E_0\times\mathbb{R}$. Since the cardinality of the
partitions $\mathcal{G}_k$
defined above is less than or equal to $ 2^{2^{k-1}}(2^{2k+1}+2)$, then
for $n \ge4$ a simple computation implies the minimal partition
\[
\mathcal{A}_n=\mathcal{P}\biggl( \bigcup_{3 \leq k \leq n-1} \mathcal{G}_k\biggr)
\]
has cardinality at most $\prod_{k=3}^{n-1} 2^{2^{k-1}}(2^{2k+1}+2) \le
\prod_{k=3}^{n-1} 2^{2^{k-1}}2^{2k+2} \leq2^{2^n}$. By
(\ref{Delta}) and Lemma~\ref{partitions} we have
\[
\sup_{(t,y)}\sum_{n\ge r}2^{n/2}\Delta_{\tau
}(A_{n}(t,y))\le C\biggl\{ \sup_{t} \sum_{n\ge r}2^{n/2}\Delta_{\rho
}(B_{n}(t))+ \sum_{n\geq r}2^{n/2}2^{-n}\biggr\}.\vadjust{\goodbreak}
\]
Thus (\ref{gamma-tau}) implies that $\tau$ satisfies (\ref{tauPreG})
with respect to the sequence of admissible partitions $\mathcal{A}_n$
on $E_0 \times\mathbb{R}$.
\end{pf}
\begin{pf*}{Proof of Theorem~\ref{mainresult}}
Let $Q$ denote the rational numbers.
Then, if we restrict the partitions $\mathcal{A}_n$ of $E_0 \times
\mathbb{R}$ in Lemma~\ref{tauPG} to $E_0 \times Q$, we immediately have
%
%
\begin{equation}\label{PreG1}
\lim_{r\to\infty}\sup_{(t,y) \in E_0\times Q}\sum_{n\ge
r}2^{n/2}\Delta
_{\tau}(A_{n}((t,y)))=0,
\end{equation}
and $(E_0\times Q, \tau)$ is totally bounded. Now let $\{G_{(s,x)}\dvtx
(s,x) \in E \times\mathbb{R}\}$ be a centered Gaussian process with
$\Ex(G_{(s,x)}G_{(t,y)})= \Pr(X_s \leq x, X_t \le y)$. Then, $G$ has
$L_2$~distance $\tau$, and by (\ref{PreG1}) and Theorem~\ref{talagrand-generic}, it is
uniformly continuous on $(E_0 \times Q, \tau)$. Hence if
$\{H_{(s,x)}\dvtx
(s,x) \in E_0 \times Q\}$ is a centered Gaussian process with
\[
\Ex\bigl(H_{(s,x)}H_{(t,y)}\bigr)= \Pr(X_s \leq x, X_t \le y)-\Pr(X_s\le x)\Pr
(X_t \le y),
\]
then
\[
\Ex\bigl(\bigl(H_{(s,x)} - H_{(t,y)}\bigr)^2\bigr) =\tau^2((s,x),(t,y))- \bigl(\Pr(X_s \le x)-
\Pr(X_t \le y)\bigr)^2.
\]
Hence the $L_2$ distance of $H$ is smaller than that of $G$, and
therefore Proposition~\ref{prop1} implies the process $H$ is uniformly continuous
on $(E_0 \times Q, d_H)$. By Lemma~\ref{tauestimate} the set $E_0 \times Q$ is dense in
$(E \times\mathbb{R}, \tau)$, and since
\[
d_H((s,x),(t,y)) \le\tau((s,x),(t,y)),
\]
we also have that $E_0 \times Q$ is dense in $(E \times\mathbb{R},
d_H)$. Thus the Gaussian process $\{H_{(s,x)}\dvtx(s,x) \in E \times
\mathbb{R}\}$ has a uniformly continuous version, which we also denote
by $H$, and since $(E, d_H)$ is totally bounded, the sample functions
are bounded on $E$ with probability one.

If $\mathcal{F}= \{ I_{X_s \le x}\dvtx(s,x) \in E \times\mathbb{R})\}$,
then since
\[
d_H((s,x),(t,y))= \rho_{P}(I_{X_s \le x}, I_{X_t \le y}),
\]
the continuity of $H$ on $(E \times\mathbb{R},d_H)$ implies condition
(ii) in Theorem~\ref{agoz} is satisfied. Since $I_{X_{t}\le y}$ is
bounded, condition (i) in Theorem~\ref{agoz} is also satisfied. Therefore,
Theorem~\ref{mainresult} follows once we verify condition (iii) of Theorem~\ref{agoz}.

To verify (iii) we use Lemma~\ref{lem4}. As before, we identify the function
$f=I_{X_s \le x} \in\mathcal{F}$ with the point $(s,x) \in E\times
\mathbb{R}$. Hence, for the centered Gaussian process
\[
\{G_f\dvtx f \in\mathcal{F}\}
\]
in (iii) of Theorem~\ref{agoz}, for $(s,x) \in E \times\mathbb{R}$, we
take the process
\[
\tilde G_{(s,x)} = G_{(s,x)} + \tilde H_{s}. 
\]
In our definition of $\tilde G$ we are assuming:

(a) $\{\tilde H_s\dvtx s \in E\}$ is a Gaussian process whose law is that
of the process $\{H_s\dvtx s \in E\}$ given in the theorem, and independent
of everything in our empirical model, and\vadjust{\goodbreak}

(b) $\{G_{(s,x)}\dvtx(s,x) \in E\times\mathbb{R}\}$ is a uniformly
continuous and sample bounded version of the Gaussian process, also
denoted by $G_{(s,x)}$, but defined above on $E_0 \times Q$. The
extension to all of $E \times\mathbb{R}$ again follows by the fact
that $E_0 \times Q$ is dense in $(E\times\mathbb{R},\tau)$.

Therefore, $\tilde G$ is sample bounded and uniformly continuous on $E
\times\mathbb{R}$ with respect to its $L_2$ distance
\[
d_{\tilde G}((s,x),(t,y))= \{ \tau^2((s,x),(t,y)) + \rho^2(s,t)\}^{1/2}.
\]
Condition (iii) of Theorem~\ref{agoz} now follows easily from Lemma~\ref{lem4}
since for $(t,y)$ fixed,
\[
\{(s,x)\dvtx\lambda((s,x),(t,y)) \le\varepsilon\}
\supseteq\{(s,x)\dvtx
d_{\tilde G}((s,x),(t,y)) \le\varepsilon\},
\]
and for a random variable $Z$ bounded by one, we have
\[
\sup_{t >0}t^2\Pr(|Z|>t)\le\Pr(|Z|>0).
\]
\upqed
\end{pf*}
%

\section{Another theorem and some examples}\label{sec7}

Let $\{X_t\dvtx t \in E\}$ be a sample continuous process such that:

\begin{longlist}[(III)]
\item[(I)]
$\sup_{t \in E}|F_t(x) - F_t(y)| \le k|x-y|^{\beta}$ for all $x,y
\in\mathbb{R}$ and some $k < \infty$ and some $\beta\in(0,1]$. Note
that this condition implies that for every $t, F_{t}$ is continuous and
hence that $\tilde{F}_{t}=F_{t}$.

\item[(III)] $|X_t - X_s| \le\Gamma\phi(s,t)$ for all $s,t \in E$, and for
some $\eta>0$ and all $x \ge x_0$
\[
\Pr(\Gamma\ge x) \le x^{-\eta}.
\]

\item[(III)] For $\beta$ as in (I), and $\eta$ as in (II), there exists
$\alpha\in(0,\beta/2)$ such that
\[
\eta\biggl(\frac{1}{\alpha}- \frac{2}{\beta}\biggr)\geq2
\]
and $(\phi(s,t))^{\alpha} \le\rho(s,t), s,t \in E$, where $\rho(s,t)$
is the $L_2$ distance of a sample bounded, uniformly continuous,
centered Gaussian process on $(E, \rho)$, which we denote by $\{
H(t)\dvtx
t\in E\}$.
\end{longlist}
%
%
\begin{theorem}\label{theor5} Let $\{X_t\dvtx t \in E\}$ be a sample continuous process
satisfying \textup{(I)--(III)} above. Then
\[
\{I_{X_{s}\le x}-\Pr(X_{s}\le x)\dvtx s\in E, x\in\mathbb{R}\}
\]
satisfies the central limit theorem in $\ell_{\infty}(E\times\mathbb
{R})$. This CLT also holds under \mbox{\textup{(I)--(II)}}, provided
\textup{(II)} is strengthened to hold for all $\eta>0$ and $x \ge
x_{\eta}$, and for some $ \alpha \in(0, \beta/2)$, we have
$(\phi(s,t))^{\alpha} \le\rho(s,t), s,t \in E$, where $\rho(s,t)$ is as
in~\textup{(III)}.
\end{theorem}
%
%
\begin{rem}\label{rem8}
If the process $\{X_t\dvtx t \in E\}$ in Theorem~\ref{theor5} is a Gaussian process,
then the CLT of Theorem~\ref{theor5} holds provided (I) is satisfied, (II) is such
that $|X_t - X_s| \le\Gamma\phi(s,t)$ for all $s,t \in E$ and
$\Gamma
< \infty$ and and for some $ \alpha\in(0, \beta/2)$, we have $(\phi
(s,t))^{\alpha} \le\rho(s,t), s,t \in E$, where $\rho(s,t)$ is as
in (III).\vadjust{\goodbreak}
\end{rem}

\begin{pf*}{Proof of Theorem~\ref{theor5}}
The theorem follows by verifying the $L$ condition in Theorem~\ref{mainresult} with
respect to $\rho$ and $\{H(t)\dvtx t \in E\}$ as given in (III). Since (I)
implies the distribution functions $F_t$ are all continuous, the
distributional transforms in (\ref{Lcondition}) are simply the
distributions themselves. Therefore,
applying (I), with $\alpha$ and $ \rho$ as given in (III), we have for
all $t \in E$ that
%
%
\begin{eqnarray}
&&{\Pr}^*\Bigl({\sup_{\{s\dvtx\rho(s,t) \le\varepsilon\}}}|F_t(X_s)- F_t(X_t)|
\ge
\varepsilon^2\Bigr)\nonumber\\[-8pt]\\[-8pt]
&&\qquad\le{\Pr}^*\biggl({\sup_{\{s\dvtx\rho(s,t) \le\varepsilon\}}}|X_s- X_t|
\ge\biggl(\frac{\varepsilon^2}{k}\biggr)^{{1}/{\beta}}\biggr).\nonumber
\end{eqnarray}
Hence (II) implies
%
%
\begin{eqnarray}
&&{\Pr}^{*}\Bigl({\sup_{\{s\dvtx\rho(s,t) \le\varepsilon\}}}|F_t(X_s)- F_t(X_t)|
\ge
\varepsilon^2\Bigr)\nonumber\\[-8pt]\\[-8pt]
&&\qquad\le\Pr\biggl( \Gamma\ge\biggl(\frac{\varepsilon^2}{k}\biggr)^{
{1}/{\beta}}
\Bigl(\sup_{\{s\dvtx\rho(s,t) \le\varepsilon\}}
\phi(s,t)\Bigr)^{-1}\biggr)\nonumber
\end{eqnarray}
and since $\alpha>0$ is such that $\eta(\frac{1}{\alpha} -\frac
{2}{\beta}) \ge2$ and
\[
\Bigl(\sup_{\{s\dvtx\rho(s,t) \le\varepsilon\}} \phi(s,t)\Bigr)^{-1} \ge\Bigl(\sup_{\{
s\dvtx\rho
(s,t) \le\varepsilon\}} \rho(s,t)\Bigr)^{{-1}/{\alpha}} \ge\varepsilon
^{{-1}/{\alpha}},
\]
(III) therefore implies that
%
%
\begin{equation}
\sup_{t \in E}{\Pr}^*\Bigl({\sup_{\{s\dvtx\rho(s,t) \le\varepsilon\}}}|F_t(X_s)-
F_t(X_t)| \ge\varepsilon^2\Bigr)\le k^{{\eta}/{\beta}} \varepsilon
^{{-2\eta}/{\beta}} \varepsilon^{{\eta}/{\alpha}} \le k^{
{\eta}/{\beta
}} \varepsilon^{2} ,\hspace*{-26pt}
\end{equation}
provided $0< \varepsilon< \varepsilon_0$ is sufficiently small to imply
$k^{{-1}/{\beta}}\varepsilon^{{2}/{\beta} - {1}/{\alpha
}} >
x_0$. To obtain the final conclusion of the theorem assume $\alpha\in
(0,\beta/2)$ and $\eta$ is sufficiently large\vspace*{1pt} that $\eta(1/\alpha
-2/\beta)>2$. Then, for $0< \varepsilon< \varepsilon_{\eta}$ sufficiently
small that $k^{{-1}/{\beta}}\varepsilon^{{2}/{\beta} -
{1}/{\alpha}} > x_{\eta}$
we again have (24). Since these estimates are uniform in $\varepsilon\in
(0,\varepsilon_ 0 \wedge\varepsilon_{\eta})$, (24) then implies the
$L$~condition, and the proof is complete.
\end{pf*}
%
%
\begin{cor}\label{cor2} Let $\{Y_t\dvtx t \in[0,T]\}$ be a sample continuous $\gamma
$-fractional Brownian motion for $0<\gamma< 1$ such that $Y_0=0$ with
probability one, and set $X_t = Y_t + Z$, where $Z$ is independent of
$\{Y_t\dvtx t \in[0,T]\}$ and has a bounded density function. Then,
\[
\{I_{X_{s}\le x}-\Pr(X_{s}\le x)\dvtx s\in[0,T], x\in\mathbb{R}\}
\]
satisfies the central limit theorem in $\ell_{\infty}([0,T]\times
\mathbb{R})$.
\end{cor}
%
%
\begin{rem}
The addition of the random variable $Z$ in the previous corollary
implies the densities of $Y_t+Z, t \in E$ are all bounded by the same\vadjust{\goodbreak}
bound as that of the density of $Z$, and hence condition (I) holds with
$\beta=1$. In particular, $Z$ is not used in any other way.
Furthermore, below we will see that something of this sort is
necessary, since we will show that the CLT of the previous corollary
fails for the fractional Brownian motion process $Y$ itself, that is,
when $Z=0$.
\end{rem}
\begin{pf*}{Proof of Corollary~\ref{cor2}}
The $L_2$ distance for $\{X_t\dvtx t \in[0,T]\}$ is given
by
%
%
\begin{equation}
\Ex\bigl((X_s -X_t)^2\bigr)^{{1}/{2}}=c_{\gamma}|s-t|^{\gamma},\qquad s,t \in[0,T],
\end{equation}
and without loss of generality we may assume the process to be
normalized so that $c_{\gamma}=1$. Furthermore, it is well known that
these processes are H\"older continuous on $[0,T]$; that is, for every
$\theta< \gamma$ we have
%
%
\begin{equation}
|X_t-X_s| \leq\Gamma|t-s|^{\theta},\qquad s, t \in[0,T],
\end{equation}
where
%
%
\begin{equation}
\Ex(\exp\{c\Gamma^2\}) <\infty
\end{equation}
for some $c>0$. That $\Gamma$ has exponential moments is due to the
Fernique--Landau--Shepp theorem, and hence the corollary follows as in
Remark~\ref{rem8}, provided we take the Gaussian process $H$
to be an $\alpha\theta$ fractional Brownian motion for any fixed
$\theta
<\gamma$ and any fixed $\alpha\in(0,\frac{1}{2})$ as $\beta=1$. Hence
the corollary is proven.
\end{pf*}
%
%
\begin{cor}Let $I=[0,T]$ and $\{Y_{(s,t)}\dvtx(s,t) \in I\times I\}$
be a
sample continuous Brownian sheet, that is, the centered Gaussian
process on $I \times I$ with covariance
$E(Y_{(s,t)}Y_{(u,v)})= (s\wedge u)( t \wedge v)$ such that with
probability one $Y_{(0,t)}=Y_{(s,0)}=0$ for $s,t \in I$. Also, set
$X_{(s,t)} = Y_{(s,t)} + Z$, where $Z$ is independent of
$\{Y_{(s,t)}\dvtx
(s,t) \in I \times I\}$ and has a bounded density function. Then,
\[
\bigl\{I_{X_{(s,t)}\le x}-\Pr\bigl(X_{(s,t)}\le x\bigr)\dvtx(s,t) \in I\times I, x\in
\mathbb{R}\bigr\}
\]
satisfies the central limit theorem in $\ell_{\infty}((I\times
I)\times
\mathbb{R})$.
\end{cor}
\begin{pf}
First of all observe that since $Z$ has a bounded density, and is
independent of the Brownian sheet $Y$, we have (I) holding with $\beta
=1$. Furthermore, from Theorem 1 in the paper \citet{yeh-wiener},
these processes are H\"older continuous on $I \times I $; that is, for
$(s,t), (u,v) \in I \times I$ and $0< \gamma<1/2$, we have
%
%
\begin{equation}
\bigl|X_{(s,t)}-X_{(u,v)}\bigr| \leq\Gamma\phi((s,t),(u,v)),
\end{equation}
where
\[
\phi((s,t),(u,v)) =\biggl(\biggl(\frac{u-s}{T}\biggr)^2 + \biggl(\frac{v-t}{T}\biggr)^2\biggr)^{
{\gamma}/{2}}
\]
and
%
%
\begin{equation}
\Ex(\exp\{c\Gamma^2\}) <\infty
\end{equation}
for some $c>0$. That $\Gamma$ has exponential moments is due to the
Fernique--Landau--Shepp theorem, and hence the corollary will follow as
in Remark~\ref{rem8}, provided we take the Gaussian process $H_{(s,t)}$
to be
%
%
\begin{equation}
H_{(s,t)} = Y_s + Z_t,\qquad (s,t) \in I \times I,
\end{equation}
where the processes $\{Y_s\dvtx s \in I\}$ and $\{Z_t\dvtx t \in I\}$ are
independent $\theta$-fractional Brownian motions.
To determine $\theta$ we fix $\gamma=1/4$, and normalizing the $Y_s$
and $Z_t$ processes suitably, we have
\[
\rho^2((s,t),(u,v))=\Ex\bigl(\bigl(H_{(u,v)} - H_{(s,t)}\bigr)^2\bigr)=
\biggl( \biggl|\frac
{u-s}{T}\biggr|^{2\theta} +\biggl|\frac{v-t}{T}\biggr|^{2\theta}\biggr).
\]
Hence for any $\alpha\in(0,1/2)$ and $\gamma=1/4$, we take $\theta
\in(0,\alpha/4)$, which implies that
\[
\phi^{\alpha}((s,t),(u,v))\le\rho((s,t),(u,v)).
\]
Since each such $\theta$ yields suitable fractional Brownian motion
choices for $Y$ and $Z$, the corollary is proven.
\end{pf}

\section{Examples where our CLT fails}\label{sec8}

\subsection{Fractional Brownian motions}\label{sec81}
Since the class of sets in our CLT arises using the Vapnik--Cervonenkis
class of half lines, one might think that perhaps if i.i.d. copies of
the process $\{X_t\dvtx t \in E\}$ satisfied the CLT in $C(E)$, then the
class of sets $\mathcal{C}$ of Theorem~\ref{mainresult} would satisfy the $\operatorname{CLT}(P)$.
Our first example shows this fails, even if the process $X_t$ is
Brownian motion on $[0,1]$ tied down at $t=0$. In this example the
process fails condition (I) in Theorem~\ref{theor5} since $\Pr(X_0=0)=1$. To prove
this we show the necessary condition for $\mathcal{C}$ to satisfy the
$\operatorname{CLT}(P)$ appearing in (ii)(a) of Theorem~\ref{delta-c} fails. More precisely, since
measurability is an issue here, the next lemma shows that there is a
countable subclass $\mathcal{C}_{Q}$ of sets in $\mathcal{C}$ such that
by Theorem~3 of \citet{talagrand-donsker-sets},
$\mathcal{C}_{Q} \notin \operatorname{CLT}(P)$.
Thus $\mathcal{C}$ fails the $\operatorname{CLT}(P)$, as otherwise all subclasses also
are in $\operatorname{CLT}(P)$.
%
%
\begin{lem}\label{BrownianMotEx}
Let $\mathcal{C}=\{C_{t,x}\dvtx0\leq t \leq1,-\infty<x< \infty\}$, where
$C_{t,x}=\{z \in C[0,1]\dvtx z(t) \leq x\}$, and assume $\{X_t\dvtx t
\in
[0,1]\}$ is a sample continuous Brownian motion tied down at zero.
Also, let $\mathcal{C}_{Q}$ denote the countable subclass of $\mathcal
{C}$ given by $\mathcal{C}_{Q} = \{ C_{t,y} \in\mathcal{C}\dvtx t,y \in Q
\}$, where $Q$ denotes the rational numbers. Then, for each integer $n
\geq1$, with probability one
\[
\Delta^{\mathcal{C}_{Q}} (\{B_1,\ldots,B_n\}) =2^n,
\]
where $\Delta^{\mathcal{C}_{Q}} (\{B_1,\ldots,B_n\})= \operatorname
{card}\{
C\cap
\{B_1,\ldots,B_n\}\dvtx C \in\mathcal{C}_{Q} \}$, and $B_1,\ldots,B_n$ are
independent copies of $\{X_t\dvtx t \in[0,1]\}$.
\end{lem}
\begin{pf}
Fix $k, 1 \leq k \leq n$, and integers $1\leq j_1<\cdots<j_k\leq n$.
The first thing we want to show is that with probability one there are
suitable $C_{t,x} \in\mathcal{C}_{Q}$ such that the $k$ functions $\{
B_{j_1},\ldots,B_{j_k}\}= C_{t,x} \cap\{B_1,\ldots,B_n\}$. Of course,
since the functions $\{B_{j_1},\ldots,B_{j_k}\}$ are random, the choice
of $C_{t,x}$ also may need be random, and for this we use the law of
the iterated logarithm (LIL). This will show that with probability one
all nonempty subsets of $\{B_1,\ldots,B_n\}$ are in $\Delta
^{\mathcal
{C}_{Q}} (\{B_1,\ldots,B_n\})$, and to get the empty set with
probability one is trivial, that is, the sample functions are
continuous on $[0,1]$, but the choice of $x$ in $C_{t,x}$ can be made
arbitrarily negative. Now for the details.

Let
\[
\mathbf{u}=(u_1,\ldots,u_n),
\]
where $u_{j_1}=u_{j_2}=\cdots=u_{j_k}=1$ and all other $u_j=2$. Then
$\|\mathbf{u}\| =\break(\sum_{j=1}^n u_j^2)^{1/2}=(4n-3k)^{1/2}$. Now set
$\mathbf{v}=(v_1,\ldots,v_n)$, where
\[
v_{j_1}=v_{j_2}=\cdots=v_{j_k}=\frac{1}{2(4n-3k)^{1/2}}
\]
and all other $v_j=\frac{1}{(4n-3k)^{1/2}}$. Then $\mathbf{v}=
\mathbf{u}/(2\|{\mathbf u}\|)$ and $\|{\mathbf v}\|=1/2$.

For $x>0$, let $Lx= \log_ex$ and set
\[
W(s)= \frac{ (B_1(s),\ldots, B_n(s))}{(2sLL(1/s))^{1/2}}
\]
for $0<s\leq1$. Then the multi-dimensional compact LIL implies with
probability one that
\[
{\liminf_{s \downarrow0}} \|{\mathbf v} - W(s)\|=0,
\]
and hence with probability one there are infinitely many rational
numbers $t \downarrow0$ such that
\[
C_{t,x(t)} \cap\{B_1,\ldots, B_n\} =\{B_{j_1},\ldots, B_{j_k}\},
\]
where $x(t) \in Q$ for $t \in Q$ and
\[
\biggl|x(t) -\frac{3(2tLL(1/t))^{1/2}}{4(4n-3k)^{1/2}}\biggr|
< \frac{(2tLL(1/t))^{1/2}}{16(4n-3k)^{1/2}}.
\]
Since $k$ and the set $\{j_1,\ldots,j_k\}$ were arbitrary, and with
probability one we can pick out the subset $\{B_{j_1},\ldots,B_{j_k}\}
$; the lemma follows as the intersection of $2^n$ subsets of
probability one has probability one.
\end{pf}

The failure of the CLT also holds for all sample continuous fractional
Brownian motions $\{X_H(t)\dvtx t \in[0,1]\}$ which are tied down at zero.
The proof of this again depends on the law of the iterated logarithm
for $n$ independent copies of this process\vadjust{\goodbreak} at $t=0$, which then allows
us to prove an analog of the previous lemma. The LIL result at $t=0$
for a single copy follows, for example, by Theorem 4.1 of \citet
{goodman-kuelbs-clustering}, and then one can extend that result to $n$
independent copies by classical proofs as in \citet
{kuelbs-strong-convergence}. The details of this last step are lengthy,
but at this stage are more or less routine in the subject, and hence
are omitted. Of course, the CLT for i.i.d. copies of these processes is
obvious as they are Gaussian.

\subsection{A uniform CLT example}\label{sec82}
In the previous examples, when the distribution function $F_t$ of $X_t$
jumped, the oscillation of the processes at that point caused a failure
of our CLT. Hence one other possible idea is that if the process $\{
X_t\dvtx t \in[0,1]\}$ is Lip-1 on $[0,1]$, then our CLT might hold. For
example, this is true for the Lip-1 process $X_t=tU, t \in[0,1]$,
where $U$ is uniform on $[0,1]$. Moreover, in this example the
densities of $F_t$ still are unbounded near $t=0$.

To see this let $X_{t,j}, j=1,\ldots,n, t \in[0,1]$, be i.i.d. copies
of $X_t=tU, t \in[0,1]$, and define
\[
W_n(C_{t,y})=\frac{1}{\sqrt n} \sum_{j=1}^n \bigl[I\bigl(X_{(\cdot),j} \in
C_{t,y}\bigr) -\Pr\bigl(X_{(\cdot),j} \in C_{t,y}\bigr)\bigr],
\]
where $\mathcal{C}=\{ C_{t,y}\dvtx t \in[0,1], y \in\mathbb{R}\}$ and
$C_{t,y}= \{z \in C[0,1]\dvtx z(t) \le y\}$.
Therefore, $W_n(C_{t,y})=0$ for all $y \in\mathbb{R}$ when $t=0$, and
also when $y/t \geq1$. Moreover, if we define $\mathcal{G}=\{
(-\infty
,r]\dvtx0 \le r \le1\}$, and
\begin{eqnarray*}
\phi(I_{C_{t,y}})&=&I_{(-\infty,1]} \qquad\mbox{if } y/t>1, 0 \le t \le1\mbox{,
 or }y=0 \mbox{ and } t=0,
\\
\phi(I_{C_{t,y}})&=&I_{(\infty,0]} \qquad\mbox{if } y/t \le0
\mbox{ but not } y=0 \mbox{ and } t=0
\end{eqnarray*}
and
\[
\phi(I_{C_{t,y}})=I_{(-\infty,y/t]} \qquad\mbox{if } 0< y/t<1,0< y< 1, 0 < t
< 1.
\]
Then $\phi(\mathcal{C})= \{I_{U \leq r}\dvtx0 \le r \le1\} \equiv
\mathcal
{G}$ and $\phi$ maps $L_2$ equivalence classes of $\mathcal{C}$ onto
$\mathcal{G}$ with respect to the law of $\{X_t\dvtx t \in[0,1]\}$ for
sets in $\mathcal{C}$, and the law of $U$ for sets in $\mathcal{G}$.
Now $\mathcal{G}$ satisfies the $\operatorname{CLT}(\mathcal{L}(U))$ by the classical
empirical CLT [e.g., see
Theorem 16.4 of \citet{billingsley}], and since $\phi$ preserves
covariances we thus have
$W_n(C_{t,y})$ converges weakly to
the Gaussian centered process
$W(C_{t,y})= Y(\phi(C_{t,y}))$ on $\mathcal{C}$, where $Y((-\infty
,s])=B(s)-sB(1)$ is the tied down Wiener process on $[0,1]$; that is,
$B(\cdot)$ is a Brownian motion.

\subsection{A Lip-1 example without the CLT}\label{sec83}
In this example we see that the Lip-1 property for $\{X_t\dvtx t \in
[0,1]\}
$ is not always sufficient for our CLT.
Here
$X_0=0$, and for $0<t\leq1$, we define
\[
X_t=t \sum_{j=1}^{\infty} \bigl(\alpha_j(t)+2\bigr)I_{E_j}(t)U,
\]
where:

\begin{longlist}
\item
$E_j= (2^{-j},2^{-(j-1)}]$ for $j \geq1$.

\item $\{\alpha_j(t)\dvtx j \geq1\}$ are independent random processes with
$\alpha_j(\cdot)$ defined on $E_j$
such that for $j \ge1$,
\[
\Pr\bigl(\alpha_j (t) = \sin(2 \pi2^{j}t), t \in E_j\bigr)=1/2
\]
and
\[
\Pr\bigl(\alpha_j (t) = \sin(2 \pi2^{j+1}t), t \in E_j\bigr)=1/2.
\]

\item $U$ is a uniform random variable on $[3/2,2]$,
independent of the $\{\alpha_j\}$.
\end{longlist}

Since the $\alpha_j$'s are zero at endpoints of the $E_j$ and we have
set $X(0)= 0$,
it is easy to see $X(t)$ has continuous paths on $[0,1]$. Moreover,
$X(t)$ is Lip-1 on $[0,1]$, and $X(t)$ has a density for each $t \in
(0,1]$, but our CLT fails.

The failure of the empirical CLT can be shown by verifying a lemma of
the sort we have above for Brownian motion, and again we see the lack
of uniformly bounded densities is a determining factor.

For each integer $n \geq1$, let $X_1,\ldots,X_n$ be independent copies
of $X$, and again take
$\mathcal{C}=\{C_{t,x}\dvtx0\leq t \leq1, \infty<x< \infty\}$, where
$C_{t,x}=\{z \in C[0,1]\dvtx z(t) \leq x\}$. Also, define $\mathcal
{C}_{Q}$ as in Lemma 
\ref{BrownianMotEx}. Then, we have the following lemma, and combined with the argument in
Section~\ref{sec81}, we see the empirical CLT fails when this $X(\cdot
)$ is used.
%
%
\begin{lem}\label{lem8} For each integer $n \geq1$, with probability one
\[
\Delta^{\mathcal{C}_{Q}} (\{X_1,\ldots,X_n\}) =2^n,
\]
where $\Delta^{\mathcal{C}_{Q}} (\{X_1,\ldots,X_n\})= \operatorname
{card}\{
C\cap
\{X_1,\ldots,X_n\}\dvtx C \in\mathcal{C}_{Q}\}$.
\end{lem}
\begin{pf} As in the proof of Lemma~\ref{BrownianMotEx}, assume one wants the $k$
functions $\{X_{i_1},\ldots,X_{i_k}\}$ with probability one, where
$1\leq i_1< i_2<\cdots<i_k\leq n$.
If we write
\[
X_i(t)= t \sum_{j=1}^{\infty}
\bigl(\alpha_{i,j}(t)+2\bigr)I_{E_j}(t)U_i,
\]
where the $\{\alpha_{i,j}\dvtx j \geq1\}$ are independent copies of $\{
\alpha_j(t)\dvtx j \geq1\}$ and $\{U_i\dvtx i \geq1\}$ are independent copies
of $U$, independent of all the $\alpha_{i,j}$'s, then
this can be arranged by taking
\[
\alpha_{i,j}(t) = \sin(2 \pi2^{j+1}t),\qquad i \in\{i_1,\ldots, i_k\},
\]
and
\[
\alpha_{i,j}(t) = \sin(2 \pi2^{j}t),\qquad i \in\{1,\ldots,n\}\cap\{
i_1,\ldots, i_k\}^c,
\]
provided we set $t=t_j= 2^{-j} + \frac{1}{4}(2^{-(j-1)} - 2^{-j})$.
The probability of this configuration on the interval $E_j$ is $1/2^n$,
and hence with probability one the Borel--Cantelli lemma implies there
are infinitely many (random) $\{t_j\downarrow0\}$ such that
\[
\alpha_{i_1,j}(t_j)=\cdots=\alpha_{i_k,j}(t_j)=0\vadjust{\goodbreak}
\]
and
\[
\alpha_{i,j}(t_j)=1
\]
for all other $i \in\{1,\ldots,n\}$.
Thus with probability one there are infinitely many rational numbers $t
\downarrow0$ such that
\[
C_{t,x(t)} \cap\{X_1,\ldots, X_n\} =\{X_{i_1},\ldots, X_{i_k}\},
\]
provided
\[
x=x(t)=\frac{17t}{4}.
\]
Of course, $x(t)$ is then also in $Q$, and since $k$ and the set $\{
i_1,\ldots,i_k\}$ were arbitrary, and we can pick out $\{
X_{i_1},\ldots
,X_{i_k}\}$ with probability one using sets in $\mathcal{C}_{Q}$, the
lemma follows as the intersection of $2^n$ subsets of probability one
has probability one.

To see $X(t)$ is Lip-1 on $[0,1]$, observe that the intervals $\{
E_j\dvtx
j\geq1\}$
are disjoint, $X(t)$ is differentiable on their interiors and an easy
computation implies
\[
{\sup_{j\geq1} \sup_{2^{-j}<t<2^{-(j-1)}}} |X'(t)| \leq[ 4 \pi+ 3]U.
\]
Furthermore, $X(t)$ is continuous on $[0,1]$, so the mean value theorem
and an elementary argument shows $X(t)$ is Lip-1 on $[0,1]$, with
Lipschitz constant bounded by $(4\pi+3)U$ with probability one.
Furthermore, since $U$ is uniform on $[3/2,2]$ and independent of the
$\{\alpha_j\}$, then $X(t)$ has a density for each $t \in(0,1]$.
\end{pf}

\subsection{Variations of the $L$ condition and the CLT}\label{sec84}

Here we produce examples where the sets $\mathcal{C}$, or more
precisely the class of indicator functions given by $\mathcal{C}$, are
$P$-pre-Gaussian, and yet $\mathcal{C} \notin \operatorname{CLT}(P)$. More importantly,
they also satisfy the modified $L$ condition, that is, we say $\{X_s\dvtx s
\in E\}$ satisfies the modified $L$ condition if for all 
$\varepsilon>0$, there exists $L< \infty$ such that
%
%
\begin{equation}\label{modLcondition}
\sup_{t \in E}\Pr^{*}\Bigl({\sup_{\rho
(s,t)\le\varepsilon}}|{F}_{t}(X_{s}) -{F}_{t}(X_{t})|>\varepsilon^{2}\Bigr)\le
L\varepsilon^{2}.
\end{equation}
Of course, if the distribution functions $\{F_t\dvtx t \in E\}$ are all
continuous, this is the $L$~condition.
Hence these examples also provide motivation for the use of the
distributional transforms in the $L$ condition of (\ref{Lcondition}) used
in Theorem~\ref{mainresult}.

Notation for the examples in this subsection is as follows.
Let $E=\{1,2,\break3,\ldots\}$, and assume $D(E) =\{z\dvtx z(t) =0 \mbox{ or
} 1,
t \in E\}$ with $\mathcal{C}=\{C_{t,y}\dvtx t \in E, y \in\mathbb{R}\}$,
where $C_{t,y}=\{z \in D(E)\dvtx z(t)\le y\}$. Then, since the
functions in
$D(E)$ take only the values zero and one, we have
\[
\mathcal{C}= \mathcal{C}_0 \cup\{ \{D(E)\}\},
\]
where $ \mathcal{C}_0=\{\tilde C_{t,0}\dvtx t \in E\}$ and for $t \in E,
\tilde C_{t,0}=\{z \in D(E)\dvtx z(t)=0\}$.
Also, let $\Sigma$ denote the minimal sigma-field of subsets of $D(E)$
containing $\mathcal{C}$, and let $P$ denote the probability on $(D(E),
\Sigma)$ such that $\Pr(\tilde C_{t,0}) =p_t$ and the events $\{
\tilde
C_{t,0}\dvtx t \in E\}$ are independent events, that is, $P$ is a product
measure on the coordinate spaces of $D(E)=\{0,1\}^{\mathbb{N}}$ with
the $t$th coordinate of $D(E)$ having the two point probability that
puts mass $p_t$ on zero, and $1-p_t$ on one.
%
%
\begin{prop}
Let $\mathcal{C}$ be defined as above. Then:

\begin{longlist}
\item
$\mathcal{C}$ is $P$-pre-Gaussian whenever $p_t= o((
\log(t+2))^{-1}) \mbox{ as } t \rightarrow\infty$.

\item $\mathcal{C} \in \operatorname{CLT}(P)$ if and only if for some $r>0$,
\[
\sum_{t=1}^{\infty} \bigl(p_t(1-p_t)\bigr)^r < \infty.
\]

\item If $p_t= (\log(t+2))^{-2}$, and $\{H(t)\dvtx t \in E\}$ consists
of centered independent Gaussian random variables with $\Ex(H(t)^2)=
(\log(t+2))^{-{3}/{2}}$, then $\{X_s\dvtx s \in E\}$ satisfies
the modified $L$ condition,
$\mathcal{C}$ is $P$-pre-Gaussian and $ \mathcal{C} \notin \operatorname{CLT}$. In
particular, in view of Theorem~\ref{mainresult}, it does not satisfy the
$L$ condition.
\end{longlist}
\end{prop}
\begin{pf}Since $\mathcal{C}$ differs from $\mathcal{C}_0$ by the
single set $D(E)$ and $P(D(E))=1$, it is easy to see that $\mathcal{C}
\in \operatorname{CLT}(P)$ if and only if $\mathcal{C}_0 \in \operatorname{CLT}(P)$.
Therefore, since the events of $\mathcal{C}_0$ are independent, Theorem
3.9.1 in Dudley [(\citeyear{Dudley-unif-clt}), page 122] implies that $\mathcal{C}_0
\in \operatorname{CLT}(P)$ if and only if for some $r>0$,
\[
\sum_{t=1}^{\infty} \bigl(p_t(1-p_t)\bigr)^r < \infty.
\]
Hence (ii) holds.

Now the centered Gaussian process $\{G_{P}(C)\dvtx C \in\mathcal{C}_0\}=
\{
G_{P}(C_{t,0})\dvtx t\in E\}$, and since the random variables $\{
G_P(C_{t,0})\dvtx t\in E\}$ are mean zero and
$\Ex(G_P(C_{t,0})^2) = p_t(1-p_t)$, we have $\mathcal{C}_0$ is
$P$-pre-Gaussian provided
\[
p_t=o\bigl(\bigl(\log(t+2)\bigr)^{-1}\bigr) \qquad\mbox{as } t \rightarrow\infty.
\]
Hence\vspace*{1pt} $\mathcal{C}$ is $P$-pre-Gaussian whenever $p_t=o((
\log(t+2))^{-1})$ as $t \rightarrow\infty$, and (i) holds.
To verify (iii) we take $p_t= (\log(t+2))^{-2}$, and $\{H(t)\dvtx
t\in E\}$ to be centered
independent Gaussian random variables with $\Ex(H(t)^2)= (\log
(t+2))^{-{3}/{2}}$. If $ \rho^2(s,t)=\Ex((H(s)-H(t))^2)$, then, for
$s \neq t$,
\begin{eqnarray*}
\rho^2(s,t)&=& \bigl(\log(s+2)\bigr)^{-{3}/{2}}+ \bigl(\log
(t+2)\bigr)^{-
{3}/{2}} \\
&\ge&\max\bigl\{ \bigl(\log(s+2)\bigr)^{-{3}/{2}}, \bigl(\log
(t+2)\bigr)^{-{3}/{2}}\bigr\}.
\end{eqnarray*}
In addition, we have
$|F_t(X_s)-F_t(X_t)|= 0$ if $X_s=X_t=0$ or $X_s=X_t=1$, and
$|F_t(X_s)-F_t(X_t)|= p_t$ if $X_t \neq X_s$.\vadjust{\goodbreak}

Therefore, for all $t \in E$ fixed and $\varepsilon>0$, we have
\[
{\mathrm{Pr}^{*}} \Bigl({\sup_{\{s\dvtx\rho(s,t)\le\varepsilon\}}}
|F_t(X_s)-F_t(X_t)|>\varepsilon^2\Bigr)= 0
\qquad\mbox{if } p_t \leq\varepsilon^2
\]
and
\begin{eqnarray*}
&&{\mathrm{Pr}^{*}} \Bigl({\sup_{\{s\dvtx\rho(s,t)\le\varepsilon\}}}
|F_t(X_s)-F_t(X_t)|>\varepsilon^2\Bigr) \\
&&\qquad\leq\mathrm{Pr}^{*}\Bigl(\sup_{\{s\dvtx
\rho
(s,t)\le\varepsilon\}} I_{X_t \neq X_s} >0 \Bigr) \qquad\mbox{for } p_t >
\varepsilon^2.
\end{eqnarray*}
Of course, $E$ countable makes the outer probabilities in the above,
ordinary probabilities, but for simplicity we retained the outer
probability notation used in (\ref{weakLcondition}) and (\ref
{Lcondition}). Now $\rho(s,t) \leq\varepsilon$ implies
\[
\max\bigl\{ \bigl(\log(s+2)\bigr)^{-{3}/{4}}, \bigl(\log
(t+2)\bigr)^{-
{3}/{4}}\bigr\} \le\varepsilon.
\]
Thus if $ p_t = (\log(t+2))^{-2}> \varepsilon^2 $, we have $\{\sup
_{\{s\dvtx\rho(s,t)\le\varepsilon\}} I_{X_t \neq X_s} >0 \} =
\varnothing$.
Combining the above we have for each fixed $t \in E$ and $\varepsilon
>0$ that
\[
{\mathrm{Pr}^{*}} \Bigl({\sup_{\{s\dvtx\rho(s,t)\le\varepsilon\}}}
|F_t(X_s)-F_t(X_t)|>\varepsilon^2\Bigr)= 0,
\]
and hence the modified $L$ condition for $\{X_s\dvtx s \in E\}$ holds. Thus
(iii) follows.
\end{pf}

%
%
\begin{appendix}\label{app}
\section*{Appendix: Talagrand's continuity result for Gaussian~processes}

The proof of Theorem~\ref{talagrand-generic} in Section~\ref{sec5} is as follows.
\begin{pf*}{Proof of Theorem~\ref{talagrand-generic}}
First we will show (i) and (ii) are equivalent. If (ii)
holds, then
by Fatou's lemma we have
\[
0 = \lim_{n \rightarrow\infty} \Ex\Bigl({\sup_{d(s,t) \le1/n}}|X_s - X_t|\Bigr)
\ge\Ex\Bigl( {\liminf_{n \rightarrow\infty} \sup_{d(s,t) \le1/n}}|X_s - X_t|\Bigr).
\]
Thus, with probability one
\[
{\liminf_{n \rightarrow\infty} \sup_{d(s,t) \le1/n}}|X_s - X_t|=0,
\]
and since the random variables ${ \sup_{d(s,t) \le1/n}}|X_s - X_t|$
decrease as $n$ increases, this implies with probability one
\[
{\lim_{n \rightarrow\infty} \sup_{d(s,t) \le1/n}}|X_s - X_t|=0,
\]
which implies (i).\vadjust{\goodbreak}

If we assume (i), then since $(T,d_X)$ is assumed totally bounded, we have
\[
Z={\sup_{t \in T}} |X_t| < \infty
\]
with probability one, and the Fernique--Landau--Shepp theorem implies
$Z$ is integrable. Since
\[
{\sup_{d(s,t) \le\varepsilon}}|X_s- X_t| \leq2Z,
\]
and (i) implies
\[
{\lim_{\varepsilon\rightarrow0} \sup_{d(s,t)\le\varepsilon}}|X_s- X_t| =0
\]
with probability one, the dominated convergence theorem implies (ii).
Thus (i) and (ii) are equivalent.

Now we assume (i) and (ii), and
choose $\varepsilon_{k}\downarrow0$ such that
\[
\sup_{s}\Ex\sup_{\{t\dvtx d(t,s)\le\varepsilon_{k}\}}X_{t} \le\sup
_{s}\Ex\sup
_{\{t\dvtx d(t,s)\le\varepsilon_{k}\}}(X_{t}-X_{s})\le\Ex\sup_{d(s,t)\le
\varepsilon_{k}}(X_{t}-X_{s}) \le2^{-k}.
\]
Since we are assuming (i) and that $(T,d)$ is totally bounded, the
sample paths of $\{X_t\dvtx t \in E\}$ are uniformly continuous and bounded
on $(T,d)$.
Hence by Sudakov's inequality, if $N(T,d,\varepsilon) $ equals the minimal
number of open balls of radius $\varepsilon$ that cover $(T,d)$, then
\[
\lim_{\varepsilon\downarrow0} \varepsilon(\log N(T,d, \varepsilon
))^{{1}/{2}}=0.
\]
Therefore, we also are free to assume the $\varepsilon_{k}$ are such that
for all $k \geq1$ we have $\varepsilon_k(\log N(\varepsilon_k))^{
{1}/{2}} \equiv\varepsilon_k( \log(N(T,d,\varepsilon_k))^{{1}/{2}} <
\frac{1}{2}$. Moreover, since $(T,d)$ is totally bounded, and (i)
holds, by Theorem 2.1.1 of \citet{talagrand-generic} there exists an
admissible sequence of partitions $\{\mathcal{\tilde A}_n\dvtx n \ge0\}$
of $(T,d)$ such that for a universal constant $L$ we have
\[
\frac{1}{2L} \sup_{t \in T} \sum_{n \geq k} 2^{n/2}\Delta(A_n(t))
\le
\Ex\Bigl(\sup_{t \in T} X_t\Bigr).
\]

Now choose $\{n_{k}\dvtx k \ge1\}$ to be a strictly increasing sequence of
integers such that $ n_1>4$ and
%
%
\begin{equation}
2^{\sum_{2 \le j\le k}1/\varepsilon_{j}^{2}}\le n_{k-1}.
\end{equation}

Based on the $n_k$'s we define an increasing sequence of partitions,
$\mathcal{B}_{n}$. For $0 \le n\le n_{1}$ we let $\mathcal
{B}_{n}=\mathcal{\tilde A}_{n}$.
For $n_{1}<n\le n_{2}$ we proceed as follows.

First we choose a maximal set $\{s_{1},\ldots, s_{N(\varepsilon_{2})}\}$
of $(T,d)$ for which $d(s_{i}$, $s_{j})\ge\varepsilon_{2}$. Furthermore, by
our choice of $\{\varepsilon_k\dvtx k \geq1\}$
via Sudakov's inequality, we have that $N(\varepsilon_{2})\le
2^{1/\varepsilon
_{2}^{2}}$.
To define the partitions for $n_1<n \le n_2$ we next consider the
partition of $T$
formed by the sets
%
%
\begin{equation}
C_j= B(s_j,\varepsilon_2) \cap\Biggl(\bigcup_{k=1}^{j-1}B(s_k, \varepsilon
_2)\Biggr)^{c},\qquad 1\le j \le N(\varepsilon_2),
\end{equation}
and the sets $B(s,\varepsilon)$ are $\varepsilon$ balls centered at
$s$. Then
by Theorem 2.1.1 of \citet{talagrand-generic} for every integer $1 \le j
\le N(\varepsilon_2)$ there exists an admissible sequence of partitions
for $(C_j,d)$, which we denote by $\mathcal{B}_{n_1,n}^{s_j}$, such that
\begin{eqnarray*}
2^{-2}&\ge&\Ex\sup_{\{t\in C_{j}\}}X_t \ge
\frac{1}{2L} \sup_{{\{t\in C_{j}\}}}\sum_{n\ge0}2^{n/2}\Delta
(B^{s_{j}}_{n_{1},n}(t)) \\
&\ge&\frac{1}{2L} \sup_{{\{t\in C_{j}\}}}\sum_{n_{1}<n\le
n_{2}}2^{n/2}\Delta\bigl(A_{n_{1}}(t)\cap B^{s_{j}}_{n_{1},n}(t)\bigr).
\end{eqnarray*}
Since the sets $C_{j}$ form a partition of $T$, if $B_{n_{1},n}$ is one
of the sets, $B^{s_{j}}_{n_{1},n}(t)$, then
\[
2^{-2}\ge\frac{1}{2L} \sup_{t\in T}\sum_{n_{1}<n\le
n_{2}}2^{n/2}\Delta\bigl(A_{n_{1}}(t)\cap B_{n_{1},n}(t)\bigr),
\]
and we define the increasing sequence of partitions $\mathcal
{B}_{n_{1},n}$ to be all sets of the form $A_{n_{1}}(t)\cap
B_{n_{1},n}(t), $ where $t \in T$ and $B_{n_1,n} \in\mathcal
{B}^{s_{j}}_{n_{1},n}(t)$ for some $j \in[1, N(\varepsilon_2)]$.
Furthermore, since the $C_j$'s are disjoint, for $n_1<n \le n_2$ we have
%
%
\begin{equation}
\card(\mathcal{B}_{n_1,n}) \le2^{2^{n_{1}}}2^{2^{n}}N(\varepsilon
_{2})\le
2^{2^{n+1}}2^{1/\varepsilon_{2}^{2}}\le2^{2^{n+1}} n_{1}\le
2^{2^{n+1}}2^{2^{n}}=2^{2^{n+2}},\hspace*{-35pt}
\end{equation}
and for $n_1<n \le n_2$ we define $\mathcal{B}_n = \mathcal{B}_{n_1,n}$.

Iterating what we have done for $n_1<n\le n_2$, we have increasing
partitions $\mathcal{B}_{n_{k-1},n}, n_{k-1}<n \le n_k$, for which
\begin{eqnarray*}
(2L)2^{-k}&\ge&\sup_{t}\sum_{n\ge0} 2^{n/2}\Delta
(B_{n_{k-1},n})\\
&\ge&\sup_{t}\sum_{n_{k-1}<n\le n_{k}}2^{n/2}\Delta
\bigl(B_{n_{k-1},n}(t)\cap
B_{n_{k-2},n_{k-1}}(t)\cap\cdots\\
&&\hspace*{138pt}{}\cap B_{n_{1},n_{2}}(t)\cap A_{n_{1}}(t)\bigr),
\end{eqnarray*}
and for $n_{k-1}<n\le n_k$ we define $\mathcal{B}_n=\mathcal{B}_{n_{k-1},n}$.
Therefore, we now have an increasing sequence of partitions $\{\mathcal
{B}_n\dvtx n \ge0\}$ such that
\fontsize{10pt}{\baselineskip}{
\begin{eqnarray*}
(2L)\sum_{k\ge r}2^{-k}&\ge&\sum_{k\ge r}\sup
_{t}\sum
_{n_{k-1}<n\le n_{k}}2^{n/2}\Delta\Biggl(B_{n_{k-1},n}(t)\cap\Biggl( \bigcap
_{j=2}^{k-1} B_{n_{j-1},n_{j}}(t)\Biggr)\cap A_{n_{1}}(t)\Biggr)\\
\hspace*{-1pt}&\ge&\sup_{t}\sum_{k\ge r}\sum_{n_{k-1}<n\le n_{k}}2^{n/2}\Delta
\Biggl(B_{n_{k-1},n}(t)\cap\Biggl( \bigcap_{j=2}^{k-1} B_{n_{j-1},n_{j}}(t)\Biggr)\cap
A_{n_{1}}(t)\Biggr)
\end{eqnarray*}}
\normalsize
and, letting $B_n(t)$ denote the generic set of $\mathcal{B}_n$
containing $t$, we have
%
%
\begin{equation}
(2L)\sum_{k\ge r}2^{-k}\ge\sup_{t}\sum_{k\ge r}\sum_{n_{k-1}<n\le
n_{k}}2^{n/2}\Delta(B_{n}(t)).
\end{equation}
Now we count the number of elements in each partition. Since (1) holds,
the partitions $\mathcal{B}_{n_{k-1},n}^{s_j}$ are assumed admissible,
and the $C_j$'s used at the subsequent iterations are always disjoint,
we have
for $n_{k-1}<n\le n_{k}$ that
\[
\card(\mathcal{B}_{n_{k-1},n})\le2^{2^{n_{1}}}\Biggl[\prod
_{j=2}^{k-1}2^{2^{n_{j}}}N(\varepsilon_{j})\Biggr]2^{2^{n}}
N(\varepsilon_{k})
\le2^{\sum_{1 \le j\le k-1}2^{n_{j}}+2^n}n_{k-1}\le2^{2^{n+2}}.
\]

Given the increasing sequence of partitions $\{\mathcal{B}_n\dvtx n \ge
0\}
$, we now define the partitions $\mathcal{A}_n$ to be the single set
$T$ for $n=0,1$ and
$\mathcal{A}_n = \mathcal{B}_{n-2}$ for $n \ge2$. Since we have
$\card
(\mathcal{B}_n) \le2^{2^{n+2}}$ for $n\ge0$ we thus have that the
$\mathcal{A}_n$'s are admissible and using (35) above they satisfy (iii).

Now assume (iii) holds and $(T,d)$ is totally bounded.
We give a sketch of the case (iii) implies (ii).
Corollary 1.6 in \citet{fernique-compactness}
reduces our task to showing
%
%
\begin{equation}\label{first}
\lim_{\eta\to0}\sup_{t}\E\sup_{s\in B_{d}(t,\eta
)}(X_{s}-X_{t})=0
\end{equation}
and
\begin{equation}
\label{second}
\lim_{\delta\to0}\delta^{2}\log_{2} N(T,\delta)=0.
\end{equation}
In the computation below the existence of $K$ follows from Theorem
2.1.1 of \citet{talagrand-generic}. To show (\ref{first}) we estimate
\begin{eqnarray*}
\sup_{t}\E \sup_{s\in B_{d}(t,\eta
)}(X_{s}-X_{t})&=&\sup
_{t}\E\sup_{s\in B_{d}(t,\eta)}X_{s}\\
&\le& K\sup_{t}\sum_{n\ge0}2^{n/2}\Delta\bigl(A_{n}(t)\cap B_{d}(t,\eta
)\bigr)\\
&\le& K\Biggl(\sup_{t}\sum_{0\le n\le k}2^{n/2}\Delta\bigl(A_{n}(t)\cap
B_{d}(t,\eta)\bigr)\\
&&\hspace*{15pt}{}+\sup_{t}\sum_{k< n}2^{n/2}\Delta\bigl(A_{n}(t)\cap
B_{d}(t,\eta)\bigr)\Biggr)\\
&\le& K\Biggl((2\eta)C2^{k/2}+\sup_{t}\sum_{k< n}2^{n/2}\Delta
\bigl(A_{n}(t)\cap
B_{d}(t,\eta)\bigr)\Biggr)\\ 
&\le& K\Biggl((2\eta)C2^{k/2}+\sup_{t}\sum_{k< n}2^{n/2}\Delta(A_{n}(t))\Biggr),
\end{eqnarray*}
where in the third inequality $C$ is such that $ \sum_{0 \le n\le
k}2^{n/2}\le C 2^{k/2}$.
Hence,
\[
\lim_{\eta\to0}\sup_{t}\E\sup_{s\in B_{d}(t,\eta
)}(X_{s}-X_{t})\le\sup_{t}K\sum_{k< n}2^{n/2}\Delta(A_{n}(t))\qquad
\mbox{for every } k.
\]
By the hypothesis, this last quantity converges to $0$ as $k\to\infty$.

To handle (\ref{second}), by the hypothesis we can choose $k$ such that
%
%
\begin{equation}\sup_{t\in T}\sum_{n\ge k}2^{n/2}\Delta(A_{n}(t))\le
\varepsilon.
\end{equation}
Hence,
\[
2^{k/2}\sup_{t\in T}\Delta(A_{k}(t))=2^{k/2}\sup_{B\in\mathcal
{A}_{k}}\Delta(B)\le\varepsilon.
\]
For each $B\in\mathcal{A}_{k}$ pick a point, $t_{k,B}$. Let $\delta
_{k}=2\sup_{B\in\mathcal{A}_{k}}\Delta(B)$. Then, if
$B_{d}(t,\delta
)=\{s\in T\colon d(s,t)\le\delta\}\mbox{ and } \frac{2\varepsilon
}{2^{k/2}}=\delta_{k}^{\prime}$, we have
\[
B(t_{k,B})\subseteq B_{d}(t_{k,B},\delta_{k}) \subseteq
B_{d}(t_{k,B},\delta_{k}^{\prime})\qquad \mbox{for every } B\in\mathcal{A}_{k}.
\]
Since $\mathcal{A}_{k}$ is a partition, $T=\bigcup_{B\in\mathcal
{A}_{k}}B_{d}(t_{k,B},\delta_{k}^{\prime})$.
Hence,
\[
\log_{2}\mathcal{N}(T,\delta_{k}^{\prime})\le\log
_{2}(2^{2^{k}})=\biggl(\frac
{2\varepsilon}{\delta_{k}^{\prime}}\biggr)^{2}.
\]
By interpolating we get (\ref{second}).
\end{pf*}
\end{appendix}

\section*{Acknowledgments}

It is a pleasure to thank the referee for his careful reading of the
manuscript. His comments and suggestions led to a number of refinements
in the presentation.


%

\printaddresses

\end{document}